\documentclass[11pt]{report}
\usepackage{sfuthesis}
 \usepackage{amssymb,latexsym,amsmath,amsthm}
 \usepackage{mathrsfs}
\def\noi{{\noindent}}

\def\bbC{{\mathbb{C}}}

\def\bbZ{{\mathbb{Z}}}

\def\MAJA{{\mathfrak{M}}}
\def\MINA{{\mathfrak{m}}}
\def\cardN{{\mathcal{N}}}
\def\SS{{\mathfrak{S}}}
\def\SI{{\mathfrak{J}}}
\def\mmod{{\text{  (mod }}}
\newtheorem{lem}{Lemma}[section] 
\newtheorem{prop}{Proposition}[section]
\newtheorem{cor}{Corollary}[section]
\newtheorem{thm}{Theorem}
\newtheorem{thm2}{Theorem}[section]
\newtheorem{define}{Definition}[section] 
\begin{document}
\title{Small Prime Solutions to Cubic Diophantine Equations}
\author{Desmond Leung}
\qualification{B.A., Wilfrid Laurier University, 2004}
\degree{Master of Science}
\dept{Mathematics}
\submitdate{Summer 2006}
\chair{Dr. Petr Lisonek}
\approvalitem{Dr. Stephen Kwok-Kwong Choi\\ Senior Supervisor}
\approvalitem{Dr. Jason Bell\\ Supervisory Committee}
\approvalitem{Dr. Yoonjin Lee\\ Supervisory Committee}
\approvalitem{Dr. Nils Bruin\\ Internal/External Examiner} 
\beforepreface 
\prefacesection{Abstract} 
\renewcommand{\baselinestretch}{1.6} \tiny\normalsize
Let $a_1, \dots, a_9$ be non-zero integers and $n$ any integer. Suppose that 
\[ a_1+\cdots +a_9 \equiv n \pmod{2}\]
and
\[ (a_i, a_j) = 1 \quad {\text{for}} \quad 1 \le i<j\le 9. \] In this thesis we will prove that\\ (i) if each of the $a_j$'s are not all of the same sign, then the cubic equation \[ a_1p_1^3 + \cdots +a_9p_9^3 =n \] has prime solutions satisfying $p_j \ll n^{1/3} + \max\{|a_j|\}^{20+\epsilon}$; and \\ (ii) if all $a_j$ are positive and $n \gg \max\{|a_j|\}^{61+\epsilon}$, then the cubic equation is soluble in primes $p_j$.\\  This result is motivated from the 2002 result for $k = 2$ by S.K.K. Choi and J. Liu. To prove the results we will use the well-known Hardy-Littlewood Circle method, which we will outline in the thesis. Lastly, we will make a note on possible generalizations of this particular problem. \\

\noindent{\bf Keywords:} Circle method, Number Theory, Goldbach's Conjecture, \\Waring's Problem, Diophantine Equations
\prefacesection{Dedication}
\vskip1.5in
\begin{center}
To truth and beauty.
\end{center}

\prefacesection{Acknowledgments}
The author would like to thank his supervisor, Dr. Stephen Kwok-Kwong Choi, for all the invaluable support he has provided over the years. One cannot ask for a better supervisor and mentor than him. 

The author would also like to thank his committee members: Dr. Nils Bruin, Dr. Yoonjin Lee, Dr. Jason Bell, as well as the committee chair Dr. Petr Lisonek for sitting in the committee and for the suggestions on the revision of the thesis.

The author would also like to thank his former professor and friend, Dr. Edward Wang, for his mentorship, even up to this day.

Lastly, the author would like to thank his Vancity colleagues and friends: Courtney Loo, Craig Cowan, Brett Hemenway, Alana McKenzie, Terry Soo, Daniel Soo, Chris Coulter, Mark Fenwick, Sarah Prien, and Ada Lam for the studying, eating, and chillaxing sessions. Their presence and companionship made this hard knock life enduring.

\renewcommand{\baselinestretch}{1.3} \tiny\normalsize
\vspace*{\fill}
\afterpreface
\renewcommand{\baselinestretch}{1.6} \tiny\normalsize
\chapter{Introduction}
\section{The Goldbach-Waring Problem}
\subsection{Goldbach's Conjecture}

The first known instance of the Goldbach's problem appeared in 1742 in the correspondence between Goldbach and Euler. It is stated as follows: \\

\noindent {\bf Goldbach's Conjecture.} {\it Every even number $n \ge 4$ is the sum of two primes, and every odd integer $n \ge 7$ is the sum of three primes.} \\ 

The even and odd numbers of prime numbers required are respectively referred to as the binary and the ternary Goldbach's conjecture. 

Several authors took a variety of approaches to tackle the binary Goldbach's conjecture. Although none of the results have been able to prove the binary case completely, the closest to date was obtained by J. R. Chen \cite{JRChen} in 1973, using ideas from sieve methods with careful analysis and treatment of the error terms. 

\begin{thm}[Chen, 1973 \cite{JRChen}] Every sufficiently large even integer $n$ can be represented in the form $n = p + P_2$, where $P_2$ is a positive integer with at most two distinct prime factors.
 \end{thm}
In the case of the ternary Goldbach's conjecture, the first significant contribution was due to Hardy and Littlewood in 1923 \cite{HL1}\cite{HL2}, by using what is now known as the {\it Hardy-Littlewood method}, or the {\it Circle method}, which we will discuss in the next chapter. Their result depended on the Generalized Riemann Hypothesis (GRH). In 1937, I. M. Vinogradov used an improved version of the Circle method to remove the dependence on GRH. 
\begin{thm}[Vinogradov, 1937 \cite{Vino2}] Every sufficiently large odd integer is the sum of three primes. \end{thm}
It is noted that this result holds only for sufficiently large value of $n$. One needs to check up to a given $n_0$ either numerically or by other methods for the proof to hold for all $n \ge 7$. During this writing, the most current improvement is due to Liu and Wang  \cite{LiuWang}, which says the Ternary Goldbach holds for $n \ge 10^{1346}$. However, at the time of this writing, to check up to $10^{1346}$ is still far from realistic to be verified by computation.  

It is important to note that under GRH, Theorem 2 will hold for $n \ge 7$.

 \subsection{Waring's Problem}

The motivation for Waring's problem stems from Lagrange's four squares theorem in 1770, which states that every positive integer is the sum of four squares. Also in 1770, Waring proposed a generalized version of the four squares problem, which is referred to as {\it Waring's problem}. \\

\noindent{\bf Waring's Problem.} {\it For every integer $k \ge 2$, there exists an integer $s = s(k)$, which depends on $k$, such that every natural number $n$ is the sum of at most $s$ $k$-th powers of natural numbers. }\\

A question was naturally raised: for any given $k$, if such an $s$ exists, what would be the least positive integer $s$? It was conjectured to be $s \le 2^k + 1$, and was first proven by Hilbert using a very intricate combinatorial argument. Later, Hardy and Littlewood gave a much simpler proof for all powers $k$ using their newly developed circle method. 

If Waring's problem was rephrased so to allow it to be true only when {\it $n$ is sufficiently large}, it turns out we can decrease the number of terms (decrease $s$) for Waring's problem to hold. One can find the most recent results in Kumchev and Tolev's expository paper \cite{Kumchev2}. 

\subsection{The Waring-Goldbach Problem}
Given the two problems listed above, for a given $k \ge 1$ a natural question to ask would be, for sufficiently large $n$ what would be the least $s = s(k)$ such that the equation
\[ n = p_1^k + p_2^k + \cdots + p_s^k \]
holds for primes $p_1, \dots, p_s$.

This is known as the {\it Waring-Goldbach Problem} and was first settled by Vinogradov \cite{Vino1} and Hua \cite{Hua1}. Over the years various authors proposed different variants of the problem, and some were solved while others are still open. One particular question was raised by Baker \cite{Baker}, which asked: if we have a $s$ fixed integers $a_1,\dots, a_s$, for sufficiently large $n \ge n_0$ with $(n, a_1, \dots, a_s) =1$, is the equation
\[ n = a_1p_1^k + \cdots +a_sp_s^k\]
soluble in primes $p_1, \dots, p_s$?

Tsang and Liu made progress for the linear case $k=1$ with $s=3$ \cite{Tsang2} and the quadratic case $k=2$ with $s=5$ \cite{Tsang1}. Later, Choi and Liu \cite{Choi1} and Choi, Liu, and Tsang \cite{Choi4} studied the same problem but with the addition of a technical condition that requires the integers $a_1, \dots, a_s$ to be pairwise co-prime. In other words,
\[ (a_i, a_j) =1,\qquad 1 \le i < j \le s. \]
With this condition much of the numerical computations can be reduced significantly. In particular, we avoid the use of Siegel zeros and the Deuring-Heilbronn phenomenon (the interested reader may refer to \cite{Choi1}). 

Choi and Kumchev also studied the linear case $k=1$ by taking a different approach, finding mean-value estimates to Dirichlet polynomials \cite{Choi2}, which can improve the estimates for the major arc.

\section{Our Specific Problem}
\noi For any integer $n$ we will consider the cubic equation

\begin{equation} \label{eq:cubicsum}
a_1p_1^3 + a_2p_2^3 + \cdots + a_9 p_9^3 = n,
\end{equation}
where $p_j$ are prime variables and the coefficients $a_j$ are non-zero integers. Either $p_j = 2$ for some $j$ or
\begin{equation} \label{eq:necsol}
a_1 + \cdots + a_9 \equiv n \!\!\!\! \pmod{2}. \end{equation}

We also suppose 
\begin{equation} \label{eq:techcond}
(a_i, a_j) = 1, \qquad (n, a_1, \dots, a_9) = 1,
\end{equation}
and denote $D = \max\{2, |a_j|, 1 \le j \le 9 \}$. Our two main results are:

\begin{thm} \label{thm:main1}
Suppose \eqref{eq:necsol} and \eqref{eq:techcond} holds. If $a_1, \dots, a_9$ are not all of the same sign, then \eqref{eq:cubicsum} has solutions in primes $p_j$, satisfying
\[ p_j \ll n^{1/3} + D^{20+\epsilon}, \]
where the implied constant depends only on $\epsilon$. 
\end{thm}

\begin{thm} \label{thm:main2}
Suppose \eqref{eq:necsol} and \eqref{eq:techcond} hold. If $a_1, \dots, a_9$ are all positive, then \eqref{eq:cubicsum} is soluble whenever 
\[ n \gg D^{61 + \epsilon}, \]
where the implied constant depends only on $\epsilon$. 
\end{thm}


We will prove these theorems using the circle method. The idea will be elaborated in Chapter 2. We will reiterate that similar to the paper \cite{Choi1}, by imposing a stronger condition \eqref{eq:techcond} than just the natural condition \eqref{eq:necsol} we need not deal with the possible existence of the Siegel zero and so the Deuring-Heilbronn phenomenon can be avoided, thus in contrast to \cite{Tsang1} or \cite{Tsang2} we can avoid much of the heavy numerical computations. 

\section{The Hardy-Littlewood Circle Method}
{\it The following were excerpted from several sources, most notably from Heath-Brown \cite{HB1}, and from Vaughan \cite{Vaughan}.} 

The idea started when Hardy and Ramanujan were working on problems involving partition functions and sums of squares (ca. 1919). Later, Hardy and Littlewood used these same ideas to prove Waring's problem. In 1937, Vinogradov refined the arguments Hardy and Littlewood used by introducing exponential sums to replace the need to ``count points on the arc of the unit circle." 

To put more in rigorous terms, let $r(n)$ be the number of solutions to the given Diophantine equation. Consider the generating function
\[ F(\alpha) = \sum_r r(n) e(\alpha n), \]
where $e(x) := \exp(2\pi i x)$. By the Fourier coefficient formula, we have
\[ r(0) = \int_0^1 F(\alpha) d\alpha. \]
If the coefficient $r(n)$ satisfy some arithmetic conditions the behaviour of $F(n)$ will be determined by an appropriate rational approximation $a/q$ to $\alpha$, with small values of $q$ usually producing large values of $F(\alpha)$. When $\alpha$ lies in an interval $[a/q - \alpha, a/q + \alpha]$ with $q$ small, a `major arc,' one hopes to estimate $F(\alpha)$ asymptoticly, while if the corresponding $q$ is large, for the `minor arcs,' one hopes that $F(\alpha)$ will be small, at least on average. 

To use this method, one uses the circle method by means of various mean-value estimates. Inequalities from Hua and Weyl (both are in Vaughan \cite{Vaughan}) give very good estimates for the minor arc, while the major arc uses information given from zeros of Dirichlet $L$-functions to break the problem into the {\it singular series} and the {\it singular integral}. We would show that the singular series has a constant contribution, while we can derive an asymptotic formula from the singular integral, and that it dominates the minor arc.


\chapter{Applying the Circle Method to our Problem}
\section{Preamble}

We consider the prime solutions of the cubic diophantine equation \eqref{eq:cubicsum},
\[ a_1p_1^3 + a_2p_2^3 + \cdots + a_9 p_9^3 = n \]
with integers $a_1,\ldots ,a_9$ and $n$ satisfying conditions \eqref{eq:necsol}
and \eqref{eq:techcond}.

We let $N$ be a large parameter and set $M$ so that
\[ CN \le M \le N, \]
for some fixed constant $0 < C < 1$. 

Let $r(n)$ be the weighted number of prime solutions of (1.1),
\[ r(n) := \sum_{\substack{n = a_1p_1^3 + \cdots + a_9p_9^3 \\ M < |a_j|p_j^3 \le N}} (\log p_1) \cdots (\log p_9). \]
Thus, $r(n)$ is a weighted count of the number of representations of $n$ as a sum of the form \eqref{eq:cubicsum} but in a restricted range.

Our aim is to ultimately show that $r(n) > 0$, but we can do better in the sense that we can derive an asymptotic formula.\\ \\
We begin by defining a cubic exponential sum over primes,
\begin{equation}
\label{eq:Salpha}
S_j(\alpha) := \sum_{M < |a_j| p^3 \le N } (\log p)  e(a_j p^3 \alpha).
\end{equation}
Our method utilizes the orthogonality of the exponential integral
\[ \int_0^1 e(\alpha x) dx = \left\{\begin{array}{ll} 1 & {\textrm{if $\alpha = 0$,}} \\ 0 & {\textrm{if $\alpha \ne 0$,}} \end{array}\right. \]
which is an important tool to count the number of solutions to any given additive equation.

Let $S(\alpha) = S_1(\alpha)\cdots S_9(\alpha)$. Since $S_j(\alpha)$ is a trigonometric polynomial, we can calculate $r(n)$ by the Fourier coefficient formula,
\begin{align*}
\int_0^1 S(\alpha) e(-\alpha n)d\alpha 
& = \int_0^1 e(-\alpha n)\prod_{j=1}^9 \sum_{M < |a_j|p_j^3 \le N}(\log p_j) e(a_jp_j^3\alpha) d\alpha \\
& = \prod_{j=1}^9 \sum_{M < |a_j|p_j^3 \le N} (\log p_j )\int_0^1 e((a_1p_1^3 + \cdots +a_9p_9^3 - n)\alpha) d\alpha \\
& = \sum_{\substack{n = a_1p_1^3 + \cdots + a_9p_9^3 \\ M < |a_j|p_j^3 \le N}} (\log p_1) \cdots (\log p_9),
\end{align*}
and so
\begin{equation}
r(n) = \int_0^1 S(\alpha)e(-n\alpha) d \alpha.
\end{equation}
Let $\epsilon > 0$ be any fixed real number and
\begin{equation} \label{eq:PQ}
L := \log N, \quad P := (N/D)^{1/10-\epsilon}.
\end{equation}
We pick $c > 0$ such that $\displaystyle{Q = \frac{N}{PL^c}}$ satisfies $2P <Q$. Here, $D = \max\{2, |a_1|, \dots, |a_9|\}$.

We need Dirichlet's Theorem on rational approximation (as in Hardy and Wright \cite{HardyWright}).
\begin{thm2}[Dirichlet] \label{thm2:ratapprox}
Given any real number $\alpha$ and any positive integer $Q$, there exist integers $a$ and $q$ with $0 < q \le Q$ such that
\[ |q\alpha - a| \le \frac{1}{Q}. \]
\end{thm2}
Then from Theorem \ref{thm2:ratapprox}, each \mbox{$\alpha \in [1/Q, 1 + 1/Q]$} may be re-written in the form 
\begin{equation}\label{eq:ratapprox}
\alpha = a/q + \lambda, \qquad |\lambda| \le 1 / (qQ),
\end{equation}
for some integers $a$ and $q$, with $1 \le a \le q \le Q$ and $(a, q) = 1$. We denote $\MAJA(q, a)$ to be the set of $\alpha$ satisfying \eqref{eq:ratapprox} and define the major arcs $\MAJA$ and the minor arc $\MINA$ as follows:
\begin{equation} \label{eq:MAs} \MAJA := \bigcup_{q \le P} \bigcup_{\substack{a = 1 \\ (a, q) = 1}}^q \MAJA(q,a), \quad \MINA := \left[ \frac{1}{Q}, 1 + \frac{1}{Q}\right] \backslash \MAJA. \end{equation}

\begin{prop} \label{prop:disj1}The major arcs $\MAJA(q,a)$ are mutually disjoint.
\end{prop}
\begin{proof}
If $a/q \ne a'/q'$, it follows from $2P < Q$ that $q + q' \le 2P < Q$. Then
\begin{eqnarray*}
\left\vert \frac{a}{q} - \frac{a'}{q'}\right\vert = \left\vert \frac{aq' - a'q}{qq'} \right\vert \ge \frac{1}{qq'} > \left(\frac{1}{q} + \frac{1}{q'}\right)\left(\frac{1}{Q}\right). 
\end{eqnarray*}
Hence, $\MAJA(q,a)$ and $\MAJA(q', a')$ are disjoint.
\end{proof}

In view of Proposition \ref{prop:disj1}, we can split $r(n)$ into its major and minor arcs,
\begin{equation} \label{eq:rint}
r(n) = \int_{1/Q}^{1 + 1/Q} S(\alpha) e(-n\alpha) d\alpha = \int_{\MAJA} S(\alpha) e(-n\alpha) d\alpha + \int_{\MINA} S(\alpha) e(-n\alpha) d\alpha. \end{equation}
As mentioned in Chapter 2, we expect that the contribution of the major arc would be dominant and that of the minor arc is negligible.

\begin{define} A {\it Dirichlet character} mod $q$ is a complex function $\chi_q(n) : \bbZ \rightarrow \bbC$ such that
\begin{align*} \chi_q(1) & = 1 \\
\chi_q(n) & = \chi_q(n + q) \\
\chi_q(m)\chi_q(n) & = \chi_q(mn)
\end{align*}
for all $m,n$ if $(q,n)=1$, and $\chi_q(n) = 0$ if $(q,n) \ne 1$. 
\end{define}

From the definition, for any Dirichlet character $\chi$ mod $q$, we define the cubic character sum to be
\[ C_\chi(a) := \sum_{k=1}^q \chi(k) e\left(\frac{ak^3}{q}\right), \qquad C_q(a) := C_{\chi_0}(a), \]
where $\chi_0$ is the principal character modulo $q$. If $\chi_1, \dots, \chi_9$ are Dirichlet character modulo $q$, we write
\begin{equation} \label{eq:Bq} B(q, \chi_1, \dots, \chi_9) := \sum_{\substack{ k=1 \\ (k, q) = 1}}^q e\left(-\frac{kn}{q}\right) \prod_{j=1}^9 C_{\chi_j}(a_jk),
\end{equation}
and
\begin{equation}\label{eq:Aq} A(q) := \frac{B(q, \chi_0, \dots, \chi_0)}{\phi^9(q)}.\end{equation}
Note that the functions $B(q, \chi_0, \dots, \chi_0)$ and $A(q)$ depend on $a_1, \dots, a_9$ and $n$ which is fixed throughout, but for the sake of conciseness, we suppress this fact in our notation. \\

We also define $F(q, \chi_1, \dots, \chi_9)$ to be a summation similar to $B(q, \chi_1, \dots, \chi_9)$, 
\begin{equation} \label{eq:Fq} F(q, \chi_1, \dots, \chi_9) := \sum_{ k=1}^q e\left(-\frac{kn}{q}\right) \prod_{j=1}^9 C_{\chi_j}(a_jk).\end{equation}

Define 
\begin{equation}
\SS(n, x) = \sum_{q \le x} A(q),
\end{equation}
and later we will show that $\SS(n, \infty)$ exists and it is called {\it the singular series}.

Denote 
\begin{align}
\cardN(q) &:= \#\{(n_1, \dots, n_9) \in \bbZ^9 : 1 \le n_i \le q, \; a_1n_1^3 + \cdots + a_9n_9^3 = n \}, \\
N(q) &:= \#\{(n_1, \dots, n_9) \in \bbZ^9 : 1 \le n_i \le q, \; (n,q)=1, \; a_1n_1^3 + \cdots + a_9 n_9^3 \equiv n \mmod q)\}.
\end{align}

By convention, we will denote the Euler Totient by $\phi(n)$, the number of divisors by $d(n)$, and the number of distinct prime factors to be $\omega(n)$. 
$p, p_1, p_2,\dots$ be always denote prime numbers, and $c_1, c_2,\ldots $ will always denote some unspecified but computable positive constants.

\section{Treatment of the Major Arcs}
\subsection{Some Preliminary Lemmas}

For each $\chi_1, \dots, \chi_9$ (mod $q$), re-write $F(q, \chi_1, \dots, \chi_9)$ as
\begin{eqnarray*}
\sum_{ k=1}^q e\left(-\frac{kn}{q}\right) \prod_{j=1}^9 C_{\chi_j}(a_jk) = \sum_{\substack{1 \le h_j \le q \\ (h_j, q)=1 \\ j = 1, \dots, 9}} \chi_1(h_1)\cdots \chi_9(h_9) \sum_{k=1}^q e\left(\frac{k(a_1h_1^3+\cdots +a_9h_9^3-n)}{q}\right).
\end{eqnarray*}
Note that
\[ \sum_{k=1}^q e\left(\frac{k(a_1h_1^3+\cdots +a_9h_9^3-n)}{q}\right) = \left\{ \begin{array}{ll} q\ & {\textrm{if $a_1h_1^3 + \cdots + a_9 h_9^3 \equiv n \!\! \pmod{q}$}}, \\ 0\quad & {\text{otherwise.}} \end{array}\right. \]
Thus
\begin{equation} \label{eq: Fq2}
F(q; \chi_1, \dots, \chi_9) = q \sum_{(q)} \chi_1 (h_1) \cdots \chi_9 (h_9),
\end{equation}
where ${\displaystyle \sum_{(q)}}$ denotes the sum over $h_1, \dots,h_9$ satisfying $1 \le h_1, \dots, h_9 \le q, (h_j, q)=1$ and ${\displaystyle \sum_{j=1}^9 a_j h_j^3 \equiv n \!\! \pmod{q}}$. In the case where $\chi_1 = \cdots= \chi_9 = \chi_0 \!\!\pmod{q}$, we get
\begin{equation} \label{eq:FqNq}
F(q; \chi_0, \dots, \chi_0) = qN(q). 
\end{equation}
In other words, we have
\begin{equation} \label{eq:newNq} N(q) = q^{-1} \sum_{k=1}^q e\left(-\frac{nk}{q}\right) \prod_{j=1}^9 C_{\chi_0}(a_jk). \end{equation}
For any prime $p$, let $s(p) := 1 + A(p).$ Notice that
\begin{eqnarray*}
\phi(p)^{-9} \sum_{\substack{k=1 \\ k | p \\ k > 1}}^p e\left(-\frac{kn}{p}\right)\prod_{j=1}^9 C_p(a_jk) & = & \phi(p)^{-9} \prod_{j=1}^9 \bigg(\sum_{h_j=1}^{p-1} e(a_j h_j^3p)\bigg) \\
& = & \phi(p)^{-9} \prod_{j=1}^9 (p-1) \\
& = & \phi(p)^{-9}\cdot \phi(p)^9 = 1.
\end{eqnarray*}
We have an analogous result for $F(q, \chi_1, \dots, \chi_9)$,
\begin{eqnarray}
s(p) & = & 1+ \phi(p)^{-9} \sum_{\substack{ k=1 \\ (k, p) = 1}}^p e\left(-\frac{kn}{p}\right) \prod_{j=1}^9 C_p(a_jk),\\ 
& = & \phi(p)^{-9} \sum_{k=1}^p e\left(\frac{-kn}{p}\right) \prod_{j=1}^9\Bigg(\sum_{h_j=1}^{p-1} e\bigg(\frac{a_jh_j^3k}{p}\bigg)\Bigg)\nonumber \\
& = & \sum_{1 \le h_1, \dots, h_9 \le p-1} \phi(p)^{-9} \sum_{k=1}^p e\Bigg(\frac{k}{p} \bigg(\sum_{j=1}^9 a_j h_j^3 - n\bigg)\Bigg) \nonumber \\
& = &  \phi(p)^{-9}  N(p)p. 
\end{eqnarray} 

\begin{lem}
Both $A(q)$ and $F(q)$ are multiplicative functions of $q$. 
\end{lem}
\begin{proof} Let $(q_1, q_2) = 1$ and $q = q_1q_2$. We write $k = k_1q_2 + k_2q_1$. Then
\begin{eqnarray*} 
A(q_1q_2) = \sum_{(k, q_1q_2)=1}^{q_1q_2} e\left(-\frac{kn}{q_1q_2}\right) \prod_{j=1}^9 C_{q}(a_jk).
\end{eqnarray*}
When $k_1, k_2$ run over the reduced residue systems modulo $q_1$ and $q_2$ respectively, $k$ will run over the reduced residue system modulo $q$.  So
\begin{eqnarray*}
A(q_1q_2) & = & \sum_{(k_1, q_1) = 1}^{q_1} \sum_{(n_2, q_2) =1}^{q_2} e\left(-\frac{n(k_1q_2+k_2q_1)}{q_1q_2}\right)  \prod_{j=1}^9 \left(\sum_{h=1}^{q_1q_2} e \left(\frac{a_jh^3(k_1q_2 + k_2q_1)}{q_1q_2}\right)\right) \\
& = & \sum_{\substack{k_1 = 1 \\ (k_1, q_1) = 1}}^{q_1} \!\!\!e\left(-\frac{k_1n}{q_1}\right)\prod_{j=1}^9 \left( \sum_{h=1}^{q_1} e\left(\frac{a_jk_1h^3}{q_1}   \right)\right) \!\!\!\!\! \sum_{\substack{k_2 = 1 \\(k_2, q_2) = 1}}^{q_2} \!\!\!e\left(-\frac{k_2n}{q_2}\right)\prod_{j=1}^9 \left( \sum_{h=1}^{q_2} e\left(\frac{a_jk_2h^3}{q_2}   \right)\right) \\
& = & A(q_1)A(q_2).
 \end{eqnarray*}
For $F(q)$, in view of \eqref{eq:newNq}, the argument for $N(q)$ is essentially the same as the one for $A(q)$. \end{proof}


\subsection{Estimation of the Character Sums}
\begin{lem} \label{lem:chineq}
Let $\chi$ \emph{(mod} $p^\alpha$\emph{)} be any non-principal character and $\alpha \ge 0$, we have the following \\
\emph{(a)} If $\chi$ is primitive, $\alpha \ge1$ and $p|a$, then $C_\chi(a) = 0$. \\
\emph{(b)} If $\chi$ is principal \emph{(mod} $p^t$\emph{)}, $p\not\vert \  a$, and $t \ge \theta + \max\{\theta, \alpha\}$, where $\theta = 1$ if $p = 3$, and $\theta = 2$ if $p \ne 3$, then $C_{\chi\chi_0}(a) = 0$. \\
\emph{(c)} $|C_\chi (a)| \le 3(3, p)(a, p^\alpha)^{1/2} p^{\alpha/2}.$ \end{lem} 
\begin{proof}
(a) Let $a' = a/p$. Write $k = u + vp^{\alpha -1}$ for $1 \le u \le p^{\alpha-1}$ and $1 \le v \le p$. Then 
\begin{eqnarray*}
C_\chi(a) & = & \sum_{k=1}^{p^\alpha} \chi(k) e\left(\frac{a'k^3}{p^{\alpha-1}}\right) = \sum_{u=1}^{p^{\alpha-1}} \sum_{v=1}^p \chi(u + vp^{\alpha-1}) e\left(\frac{a'(u+vp^{\alpha-1})^3}{p^{\alpha-1}}\right) \\
& = & \sum_{u = 1}^{p^{\alpha-1}} e\left(\frac{a'u^3}{p^{\alpha-1}}\right) \sum_{v=1}^p \chi(u+ vp^{\alpha-1}).
\end{eqnarray*}
Since the inner-sum over $v$ is exactly zero, we have that $C_\chi(a) = 0$.\\

(b) For $1 \le h \le p^t$, write $h = u + vp^{t - \theta},$ where $1 \le u \le p^{t - \theta}, 1 \le v \le p^\theta$. Since $t - \theta \ge \max\{\theta, \alpha\}$, we have $h^3 \equiv u^3 + 3u^2vp^{t-\theta}$ (mod $p^t$) and $h \equiv u$ (mod $p^\alpha$). Then we get
\[ C_{\chi\chi_0}(a) = \sum_{u=1}^{p^{t-\theta}} \chi\chi_0 (u) e\left(\frac{au^3}{p^t}\right) \sum_{v=1}^{p^\theta} e\left(\frac{3au^2v}{p^\theta}\right). \]
If $(u, p) > 1$, then $\chi\chi_0(u) = 0$. On the other hand, if $(u, p) = 1$, then the inner-sum over $v$ will equal to zero. In either case, $C_{\chi\chi_0}(a) = 0$, as required.
\\

(c) If $\chi$ is primitive and $p|a$, then $C_\chi(a) = 0$ from part (a). If $p \nmid a$, we have
\begin{eqnarray}
\sum_{a=1}^{p^\alpha} |C_\chi(a)|^2 & = & \sum_{k,h=1}^{p^\alpha} \overline{\chi}(k)\chi(h) \sum_{a=1}^{p^\alpha} e\left(\frac{a(k^3-h^3)}{p^\alpha}\right) \nonumber\\
\label{eq:squarechi}
& = & p^\alpha \sum_{k^3 \equiv h^3 (\text{mod } p^\alpha)} \overline{\chi}(k)\chi(h).
\end{eqnarray}
If $X^3 \equiv 1$ (mod $p^\alpha$) and $p \ne 2$, then $(X-1)(X^2+X+1) \equiv 0$ (mod $p^\alpha$). 
Completing the square on the above quadratic polynomial yields
\[ (X+\overline{2})^2 - \overline{2}^2 + 1 \equiv 0 \pmod{p^\alpha}, \]
where $2\cdot \overline{2} \equiv 1 \pmod{p^\alpha}$. This is equivalent to
\[ (2X+1)^2 \equiv -3 \pmod{p^\alpha}. \]
On the other hand, if ${\displaystyle \left(\frac{-3}{p}\right)} \ne 1$, there is only one solution to $x^3 \equiv 1 \pmod{p^\alpha}$. Hence, \eqref{eq:squarechi} will give
\[ \sum_{a=1}^{p^\alpha} |C_\chi(a)|^2 = p^\alpha\sum_{k^3 \equiv h^3 {\textrm{ (mod $p^\alpha$)}}} \overline{\chi}(k)\chi(h) = p^\alpha\sum_{\substack{k=1 \\ (k, \; p^\alpha) = 1}}^{p^\alpha} \overline{\chi}(k)\chi(k) = \phi(p^\alpha) p^\alpha. \]

If ${\displaystyle \left(\frac{-3}{p}\right) = 1}$, we will get two incongruent solutions, namely
\[ X \equiv \overline{2}(1 \pm b) \!\!\!\pmod{p}, \]
where $b^2 \equiv -3 \pmod{p}$.

Denote $\alpha_1 = \overline{2}(1 + b) \pmod{p}$, and $\alpha_2 =\overline{2}(1 - b) \pmod{p}$.
If $n^3 \equiv m^3 \pmod{p}$, then $n \equiv m$, $\alpha_1m,$ or $\alpha_2 m \pmod{p}$. Thus from \eqref{eq:squarechi}, 
\[ \sum_{a=1}^p |C_\chi(a)|^2 = p (1 + \overline{\chi} (\alpha_1) + \overline{\chi} (\alpha_2) ) \sum_{\substack{n=1 \\ (n, p) = 1}}^{p} |\chi(n)|^2 \le 3 \phi(p)p. \] 

If $p \ne 3$ and $\alpha \ge 1,$ denote $f(X) := X^3 -1.$ Then we have $f'(X) = 3X^2.$ The solutions to $f'(X) = 3X^2 \equiv 0$ (mod $p$), or $X \equiv 0$ (mod $p$) is clearly different from the solutions to $X^3 -1 \equiv 0$ (mod $p$). It follows that $f'(X) \equiv 0 \pmod{p}$ and $f(X) \equiv 0 \pmod{p}$ share no common root. We need a Theorem of Hensel's (pg. 33 of Hua \cite{Hua} or in Apostol \cite{Apostol}) and another theorem regarding the number of solutions to polynomial congruences (pg. 32 of Hua \cite{Hua}).

\begin{thm2}[Hensel] \label{thm:Hensel} Let $f(X) = a_nX^n + \cdots + a_1X + a_0$ and $f'(X) = na_nX^{n-1} + \cdots + 2a_2X + a_1$. If $f(X) \equiv 0$ and $f'(X) \equiv 0$ \emph{(mod} $p$\emph{)} have no common solution, then the two congruences, {\mbox{$f(X) \equiv 0$ \emph{(mod} $p^\alpha$\emph{)}}} and $f(X) \equiv 0$ \emph{(mod} $p$\emph{)} have the same number of solutions.
\end{thm2}
\begin{thm2}
\label{thm:polycong} 
Let $p$ be a prime number. The number of solutions (including repeated ones) to the congruence 
\[ f(x) = a_n X^n + \cdots + a_0 \equiv 0 \!\!\pmod{p} \]
does not exceed $n$.
\end{thm2}
We will continue proving Lemma \ref{lem:chineq} (c).

In view of Theorem \ref{thm:Hensel}, the number of solutions to $f(X) \equiv X^3 -1 \equiv 0$ (mod $p^\alpha$) is the same as $f(X) \equiv X^3 -1 \equiv 0$ (mod $p$). 
So by Theorem \ref{thm:polycong}, the number of solutions is at most 3. 

By the same argument as above and from \eqref{eq:squarechi} we have
\[ \sum_{a=1}^{p^\alpha}\left|C_{\chi}(a)\right|^2=p^\alpha \sum_{n^3 \equiv m^3 (\text{mod } p^\alpha)} \overline{\chi}(n)\chi(m) \le 3\phi(p^\alpha) p^\alpha. \] 
Now, we know that if $p \equiv 1$ (mod 3), ie., when ${\displaystyle \left(\frac{-3}{p}\right)=1}$, there are ${\displaystyle \frac{p-1}{3}}$ cubes (that is, there are at most 3 solutions to $X^3 \equiv a$ (mod $p$)), and if $p \equiv 2$ (mod 3), ie., ${\displaystyle \left(\frac{-3}{p}\right)=-1}$, there are $p-1$ cubes (there is exactly one solution to {\mbox{$x^3 \equiv a$ (mod $p$)).}} 

The cubic residues have the form of $1^3, 2^3, \dots, (p-1)^3,$ and the two possible forms of a cubic nonresidues are $b\cdot 1^3, b \cdot 2^3, \dots, b \cdot (p-1)^3$ and $b^2 \cdot 1^3, b^2 \cdot 2^3, \dots, b^2 \cdot (p-1)^3$, and $b$ is not a cubic residue with $(b, p) = 1$. For the cubic residues, after substituting $n \mapsto \overline{a} n$ we have
\begin{eqnarray*}
C_{\chi} (a^3) = \sum_{n=1}^{p^\alpha} \chi(n) e\left(\frac{a^3n^3}{p^\alpha}\right) = \sum_{n=1}^{p^\alpha} \chi(\overline{a}n) e\left(\frac{n^3}{p^\alpha}\right) = \chi(\overline{a})C_\chi (1),
\end{eqnarray*}
and so $|C_\chi(a^3)| = |C_\chi(1)|$. 

Similarly, for the first type of cubic non-residues (ie. of the form $b a^3$), after substituting $n \mapsto \overline{a}n$ we have
\begin{eqnarray*}
C_{\chi} (ba^3) = \sum_{n=1}^{p^\alpha} \chi(n) e\left(\frac{ba^3n^3}{p^\alpha}\right) = \sum_{n=1}^{p^\alpha} \chi(\overline{a}n) e\left(\frac{bn^3}{p^\alpha}\right) = \chi(\overline{a})C_\chi (b),
\end{eqnarray*}
so we have $|C_\chi(ba^2)| = |C_\chi(b)|$ and $|C_\chi(b^2a^3)| = |C_\chi(b^2)|.$ 
Therefore, we can rewrite \eqref{eq:squarechi} into its cubic residues and nonresidues as
\begin{eqnarray} \label{eq:cubicsplit}
\sum_{a=1}^{p^\alpha} |C_\chi(a)|^2 & = & \frac{\phi(p^\alpha)}{3}(|C_\chi(a^3)|^2 + |C_\chi(ba^3)|^2 + |C_\chi(b^2a^3)|^2) \\
& \le & 3 \phi(p^\alpha) p^\alpha. \nonumber
\end{eqnarray}
For the first term, 
\[ \phi(p^\alpha)|C_\chi(a^3)|^2/3 \le 3\phi(p^\alpha)p^\alpha, \]
 \[ {\text{ie.,}}\qquad |C_\chi(a^3)|^2 \le 9 p^\alpha, \ \]
 \[ {\text{ie.,}}\qquad   |C_\chi(a^3)| \le 3 p^{\alpha/2}. \]
The other two terms can be shown similarly,
\[ |C_\chi(ba^3)| \le 3p^{\alpha/2} \quad {\text{and}} \quad |C_\chi(b^2a^3)| \le 3p^{\alpha/2}. \]
Hence, $C_\chi(a) \le 3p^{\alpha/2} = 3(a, p^\alpha)^{1/2}p^{\alpha/2}$.

Now, consider the case when $p^\alpha | a$. Write $k = a/p^\alpha$. Then
\[ |C_\chi(a)| = \bigg\vert \sum_{h=1}^{p^\alpha} \chi(h)e(kh^3)\bigg\vert = \bigg\vert \sum_{h=1}^{p^\alpha} \chi(h) \bigg\vert
= \left\{ \begin{array}{ll} 0 \quad & \textrm{if $\chi \ne \chi_0$}, \\
\phi(p^\alpha) \quad & \textrm{if $\chi = \chi_0$}.
\end{array} \right.
\]
On the other hand,
\[ 3(a, p^\alpha)^{1/2}p^{\alpha/2} = 3(kp^\alpha, p^\alpha)^{1/2} p^{\alpha/2} =3p^\alpha \ge \phi (p^\alpha ), \]
and so $|C_\chi(a)| \le 3(a, p^\alpha)^{1/2}p^{\alpha/2}$.

Lastly, consider the case when $p^\omega \Vert a$, with $0 < \omega < \alpha$, and $a = a'p^\omega$. 
For $1 \le n \le p^\alpha$, write $k = u + vp^{\alpha - \omega}$, where $0 < u \le p^{\alpha-\omega}$ and $0 < v \le p^\omega$. Then $k$ runs through the complete reduced residue system mod $p^\alpha$. So after the substitution we get
\begin{eqnarray*}
C_\chi(a) & = & \sum_{k=1}^{p^\alpha} \chi(k)e\left(\frac{ak^3}{p^\alpha}\right) = \sum_{k=1}^{p^{\alpha - \omega}} \chi(k)e\left(\frac{a'k^3}{p^{\alpha - \omega}}\right) \\
& = & \sum_{u=1}^{p^{\alpha-\omega}} \sum_{v=1}^{p^\omega}  \chi(u + vp^{\alpha - \omega}) e\left(\frac{a'(u + vp^{\alpha-\omega})^3}{p^{\alpha - \omega}}\right) \\
& = & \sum_{u=1}^{p^{\alpha - \omega}} \sum_{v=1}^{p^\omega} \chi(u + vp^{\alpha - \omega}) e\left(\frac{a'(u^3 + 3u^2vp^{\alpha - \omega} + 3uv^2p^{2(\alpha - \omega)} + v^3p^{3(\alpha - \omega)})}{p^{\alpha - \omega}}\right) \\
& = & \sum_{u=1}^{p^{\alpha-\omega}} e\left(\frac{a'u^3}{p^{\alpha-\omega}}\right) \sum_{v=1}^{p^\omega} \chi(u + vp^{\alpha - \omega}). 
\end{eqnarray*}
Let $\chi^*$ \!(mod $p^\eta$) be a primitive character which induces $\chi$. 
If $\eta > \alpha - \omega$ then
\begin{eqnarray*}
\sum_{v = 1}^{p^\omega} \chi(u + vp^{\alpha - \omega}) = \sum_{v = 1}^{p^\omega} \chi^*(u + vp^{\alpha-\omega}) = p^{\alpha - \eta} \sum_{v = 1}^{p^{\omega - \alpha + \eta}} \chi^*(u + vp^{\alpha - \omega}) = p^{\alpha - \eta} \cdot 0 = 0,
\end{eqnarray*}
and so $C_\chi (a) = 0$.

If $\eta \le \alpha - \omega$ then
\begin{eqnarray*}
\sum_{v = 1}^{p^\omega} \chi(u + v p^{\alpha - \omega}) = p^\omega \chi^*(u)\chi_0(u),
\end{eqnarray*}
and so
\begin{eqnarray*}
C_\chi(a) = \sum_{u =1}^{p^{\alpha - \omega}} e\left(\frac{a'u^3}{p^{\alpha - \omega}}\right) \sum_{v=1}^{p^\omega} \chi(u + vp^{\alpha - \omega}) = p^\omega \sum_{u=1}^{p^{\alpha - \omega}} e\left(\frac{a'u^3}{p^{\alpha - \omega}}\right) \chi^*(u)\chi_0(u) = p^\omega C_{\chi^* \chi_0} (a'),
\end{eqnarray*}
with $\chi_0$ \!\!\!\! (mod $p^{\alpha - \omega}$). 
Since $p \nmid a'$, we can use the bound proven in the first case and obtain
\begin{eqnarray*}
|C_{\chi'\chi_0}(a')| \le 3p^{(\alpha-\omega)/2}.
\end{eqnarray*}
Hence,
\begin{eqnarray*}
|C_\chi(a)| \le p^\omega \cdot 3p^{(\alpha - \omega)/2} = 3(a, p^\alpha)^{1/2}p^{\alpha/2}.
\end{eqnarray*}
This completes the proof of Lemma \ref{lem:chineq} (c).
\end{proof}

\begin{lem} \label{lem:convbdd}
Let $\chi_j$ \!\emph{(mod} $r_j$\emph{)} with $j = 1, \dots, 9$ be primitive characters, $\chi_0$ be the principal character \emph{(mod} $q$\emph{)}, and $r_0 = {\emph{\textrm{lcm}}}[r_1, \dots, r_9]$. Then
\[
\sum_{\substack{q \le x \\ r_0 | q}} \frac{1}{\phi^9(q)}|B(q, \chi_1\chi_0, \dots, \chi_9\chi_0)| \ll r_0^{-3+\epsilon}. \]
\end{lem}
\begin{proof} Lemma \ref{lem:chineq} (c) asserts that for any character $\chi$ \!(mod $p^\alpha$) with $\alpha \ge 0$, we have
\[ |C_\chi(a)| \le 3(3,p)(a, p^\alpha)^{1/2}p^{\alpha/2}. \]
Therefore, for character $\chi_1, \dots, \chi_9$ (mod $p^\alpha$) and from \eqref{eq:Aq}
\begin{equation*} 
|B(p^\alpha, \chi_1, \dots, \chi_9)| \le p^\alpha(3(3,p)p^{\alpha/2})^9 \prod_{j=1}^9 (a_j, p^\alpha)^{1/2} \le 3^{18} p^{6\alpha},
\end{equation*}
where in the last inequality we have used the condition \eqref{eq:techcond} that $(a_i,a_j)=1$; in fact, 
\[\prod_{j=1}^9 (a_j, p^\alpha)^{1/2} \le p^{\alpha/2}.\] Since $|B(q, \chi_1, \dots, \chi_9)|$ is multiplicative,
\[ |B(q, \chi_1, \dots, \chi_9)| = \prod_{p^\alpha \Vert q} |B(p^\alpha, \chi_, \dots, \chi_9)| \le q^63^{18\omega (q)} \le q^{6}d^{18}(q),
\] where we used the fact that $3^{\omega(q)} \le d(q)$ from Hardy and Wright \cite{HardyWright}. Thus, we have
\[ \sum_{\substack{q \le x \\ r_0|q}} \frac{1}{\phi^9(q)} |B(q, \chi_1\chi_0, \dots, \chi_9\chi_0)| \ll \sum_{\substack{q \le x \\ r_0 | q}} \frac{q^6d^{18}(q)}{\phi^9(q)}. \]
Since 
\[ \frac{k}{\phi(k)} = \prod_{p|k}\left(1 - \frac{1}{p}\right)^{-1} \ll \prod_{p|k} \left( 1 + \frac{1}{p}\right) \le \sum_{d|k}\frac{1}{d} \ll d(k),\]
taking $q = r_0k$ obtains
\begin{eqnarray*} \sum_{\substack{q \le x \\ r_0|q}} \frac{q^6 d^{18}(q)}{\phi^{9}(q)} 
& \ll & \sum_{\substack{q \le x \\ r_0 | q}}\frac{d^{27}(q)}{q^3} \ll \sum_{1 \le k \le x/r_0} \frac{d^{27}(kr_0)}{(kr_0)^3} \\
& \ll & \frac{d^{27}(r_0)}{r_0^3}\sum_{k \le x/r_0}\frac{d^{27}(k)}{k^3} \ll r_0^{-3+\epsilon},\end{eqnarray*}
because $d(q) \ll q^\epsilon$, for some $\epsilon > 0$.
\end{proof}
\subsection{Building the Asymptotic Formula}
\begin{lem} \label{lem:singularbdd}
Let $\MAJA$ be the major arcs defined in \eqref{eq:MAs}. Then we have
\begin{eqnarray*}
 \int_\MAJA S_1 (\alpha) \cdots S_9(\alpha) e(-n\alpha) d \alpha - \frac{1}{3^9}\SS(n, P)\SI(n) & \ll & \frac{N^{2}}{|a_1\cdots a_9|^{1/3}L^A}, \end{eqnarray*}
where $A>0$ is some constant and 
\[ \SI(n) := \sum_{\substack{a_1m_1 + \cdots + a_9m_9 = n \\ M < |a_j|m_j^3 \le N}} (m_1 \cdots m_9)^{-2/3}. \]
\end{lem} 
\begin{proof} 
For $j = 1, \dots, 9,$ set
\[ N_j := N / |a_j|, \qquad M_j := M/|a_j|, \qquad V_j(\lambda) := \!\!\!\!\! \!\! \sum_{M < |a_j|m^3 \le N} \!\!\!\!\!\!\! e(a_jm^3\lambda), \]
and 
\begin{eqnarray} W_j (\chi, \lambda) & := & \sum_{M < |a_j|p^3 \le N} \!\!\!\! (\log p) \chi(p) e(a_jp^3\lambda) - \delta_\chi \!\!\!\!\! \sum_{M < |a_j|m^3 \le N} e(a_j m^3 \lambda), \label{eq:Wvlog}
 \end{eqnarray}
where 
\[ \delta_\chi := \left\{\begin{array}{ll} 1 \ & {\textrm{if $\chi = \chi_0$,}} \\ 0 \ & {\textrm{otherwise.}} \end{array}\right. \]  Note that 
when $q \le P$ and $M < |a_j|p^3 \le N$, we have $(q, p) = 1$, and so $W_j(\chi_j \chi_0, \lambda) = W_j(\chi_j, \lambda)$ 
for primitive characters $\chi_j$.

By introducing Dirichlet characters, we can rewrite the exponential sum $S_j(\alpha)$ as 
\[ S_j\left(\frac{h}{q} + \lambda\right) = \frac{C_{\chi_0}(a_jh)}{\phi(q)}V_j(\lambda) + \frac{1}{\phi(q)}\sum_{\chi \!\!\!\!\!\! \pmod{q}} C_\chi(b_jh)W_j(\chi, \lambda) =: T_j + U_j, \]
(see for example, (2) \S26, \cite{Davenport1}) where $T_j$ and $U_j$ is to denote the first and second term, respectively. The first and second term can be re-interpreted as the explicit formula for finding zeros of the Dirichlet $L$-functions over the principal and primitive characters, respectively. Substituting into the major arc integral gives
\begin{eqnarray}  && \int_\MAJA S_1(\alpha) \cdots S_9(\alpha)e(-n\alpha)d\alpha \nonumber \\
&& = \sum_{q \le P} \sum_{\substack{h=1 \\ (h, q)=1}}^q e\left(-\frac{hn}{q}\right) \int_{\MAJA(a,q)} S_1\left(\frac{h}{q}+\lambda\right) \cdots S_9\left(\frac{h}{q}+\lambda\right)d\lambda \nonumber \\ \label{eq:split_terms}  && = I_0 + \dots + I_9,  \end{eqnarray}
where $I_\nu$ is the contribution from those products with $\nu$ pieces of $U_j$ and $9-\nu$ pieces of $T_j$, 
\begin{equation} \label{eq:Ivsum} I_\nu := \sum_{q \le P} \sum_{\substack{h=1 \\ (h, q) = 1}}^q e\left(-\frac{hn}{q}\right)\int_{-1/(qQ)}^{1/(qQ)}\sum_{\substack{I \subseteq \{1,2,\dots, 9\} \\ |I|= \nu}}\left(\prod_{i \in I} U_i\right) \Bigg(\prod_{j \in \{1,\dots, 9\} \backslash I}T_j\Bigg)   e(-n\lambda)d\lambda. \end{equation}
We expect that $T_j$ will contribute to the main term while $U_j$ will be negligible. In fact, 
we will show that $I_0$ gives the main term while $I_1, \dots, I_9$ will contribute to the error term. 

In view of \eqref{eq:Ivsum}, we can reduce the characters in $I_9$ into primitive characters, 
\begin{eqnarray*}
|I_9| & = & \left| \sum_{q \le P}\sum_{\chi_1 \text{(mod } q)} \cdots \sum_{\chi_9 \text{(mod } q)} \frac{B(q, \chi_1, \dots, \chi_9)}{\phi^9(q)} \int_{-1/(qQ)}^{1/(qQ)} W_1(\chi_1, \lambda) \cdots W_9(\chi_9, \lambda)e(-n\lambda) d\lambda \right| \\
& \le & \sum_{r_1 \le P} \cdots \sum_{r_9 \le P} {\sum_{\chi_1 \!\!\!\!\!\!\pmod{r_1}}}^{\!\!\!\!\!\!\!*} \cdots {\sum_{\chi_9 \!\!\!\!\!\!\pmod{r_9}}}^{\!\!\!\!\!\!\!*}  \ \ \sum_{\substack{q \le P \\ r_0 | q}} \frac{|B(q, \chi_1\chi_0, \dots, \chi_9 \chi_0)|}{\phi^9(q)}  \\
& & \times \int_{-1/(qQ)}^{1/(qQ)}|W_1(\chi_1\chi_0, \lambda)| \cdots |W_9(\chi_9\chi_0, \lambda)|d\lambda,
\end{eqnarray*}
where $\chi_0$ is the principal character modulo $q$ and $r_0 = $ {\textrm{lcm}}$[r_1, \dots, r_9]$. As mentioned previously, since $q \le P$ and $M < |a_j|p^3 \le N$, we have $(q, p) = 1$, and so \mbox{$W_j(\chi_j \chi_0, \lambda) = W_j(\chi_j, \lambda)$} for the primitive characters $\chi_j$ above. Consequently, by Lemma \ref{lem:convbdd} we obtain
 \begin{align*}
|I_9| & \le  \sum_{r_1 \le P} \cdots \sum_{r_9 \le P} {\sum_{\chi_1 \!\!\!\!\!\! \pmod{r_1}}}^{\!\!\!\!\!\!\!*} \cdots \!\!\!\! \! {\sum_{\chi_9 \!\!\!\!\!\!\pmod{r_9}}}^{\!\!\!\!\!\!\!*} \; \int_{-1/(r_0Q)}^{1/(r_0Q)} |W_1(\chi_1, \lambda)| \cdots |W_9(\chi_9, \lambda)| d \lambda \\
&  \times \sum_{\substack{q \le P \\ r_0 |q}} \frac{|B(q, \chi_1\chi_0, \dots, \chi_9\chi_0)|}{\phi^9(q)} \\
& \ll  \sum_{r_1 \le P} \cdots  \sum_{r_9 \le P} r_0^{-3+\epsilon} {\sum_{\chi_1\!\!\!\!\!\! \pmod{r_1}}}^{\!\!\!\!\!\!\!*} \cdots \!\!\!\! {\sum_{\chi_9 \!\!\!\!\!\! \pmod{r_9}}}^{\!\!\!\!\!\!\!*} \;\int_{-1/(r_0Q)}^{1/(r_0Q)} |W_1(\chi_1, \lambda)|\cdots |W_9(\chi_9, \lambda)| d \lambda.
\end{align*}
We have the bound
\[ r_0^{-3+\epsilon} = {\text{lcm}}[r_1, \dots, r_9]^{-3+\epsilon} \le r_1^{-1/3+\epsilon}\cdots r_9^{-1/3+\epsilon}, \]
which can be shown by elementary means.

We define the $L_2$ norm and $\sup$-norm estimates for $W_j$.
\begin{equation} \label{eq:Kjdef} K_j := \sum_{r \le P} r^{-1/3+\epsilon} {\sum_{\chi \!\!\!\!\! \pmod{r}}}^{\!\!\!\!\!*} \Bigg( \int_{-1/(rQ)}^{1/(rQ)} |W_j(\chi, \lambda)|^{2} d \lambda\Bigg)^{1/2}, \end{equation}
and
\begin{equation} \label{eq:Jjdef} J_j := \sum_{r \le P} r^{-1/3 + \epsilon} {\sum_{\chi \!\!\!\!\! \pmod{r}}}^{\!\!\!\!\!*} \max_{|\lambda| \le 1/(rQ)} |W_j(\lambda, \chi)|, \end{equation}
where ${\displaystyle {\sum}^*_{\chi \!\! \pmod{r}}}$ is over all the primitive characters modulo $r$.  We expect $K_j$ will be considerably smaller than $J_j$ due to
the cancellation in the integral.

\begin{lem} \label{lem:KJlem} \emph{(i)} For any fixed $A >0$, we have
\begin{equation}  \label{eq:Keq} K_j \ll_A |a_j|^{-1/3}N^{-1/6}L^{-A}.
\end{equation}
\emph{(ii)} For any fixed $A >0$, we have
\begin{equation}  \label{eq:Jeq} J_j \ll_A N_j^{1/3}L^{-A}.
\end{equation}
\end{lem}
\noindent The proof of Lemma \ref{lem:KJlem} will be given in section 3.2.5. 

Using these estimates and Cauchy's inequality, we obtain
\begin{align} 
|I_9|  \ll & \ K_1K_2 \prod_{j=3}^9 J_j \nonumber \\
\ll & \  \prod_{j=1}^2 \left\{ \sum_{r_j \le P} r_j^{-1/3+\epsilon} {\sum_{\chi_j \!\!\!\!\!\! \pmod{r_j}}}^{\!\!\!\!\!*} \left(\int_{-1/(r_jQ)}^{1/(r_jQ)} |W_1(\chi_j, \lambda)|^2d \lambda\right)^{1/2} \right\} \nonumber \\
& \times\prod_{j=3}^9 \left\{ \sum_{r_j \le P} r_j^{-1/3+\epsilon} {\sum_{\chi_j \!\!\!\!\!\! \pmod{r_j}}}^{\!\!\!\!\!*} \max_{|\lambda| \le 1/(r_jQ)} |W_j(\chi_j, \lambda)| \right\} \nonumber \\
\ll &\  \frac{N^{7/3}N^{-1/3}}{|a_1\cdots a_9|^{1/3}L^A} \nonumber \\ \label{eq:L9} \ll &\  \frac{N^2}{|a_1\cdots a_9|^{1/3}L^A}.
\end{align}
For the other $I_{\nu}$, we need an estimate for $V_j$.
\begin{lem} We have
\begin{equation} \label{eq:VjTrivial}
V_j(\lambda) = \sum_{M_j^{1/3} < m \le N_j^{1/3}} e(a_jm^3\lambda) \ll N_j^{1/3},
\end{equation}
and
\begin{equation} \label{eq:VjL2}
H_j := \left\{\int_{-1/Q}^{1/Q} |V_j(\lambda)|^2 d\lambda \right\}^{1/2} \ll N^{-1/6} |a_j|^{-1/3}. 
\end{equation}
\end{lem}
\begin{proof} The proof for \eqref{eq:VjTrivial} is trivial. For \eqref{eq:VjL2}, we have
\begin{eqnarray} \label{eq:computeI0}
V_j(\lambda) & = & \!\!\!\!\! \! \sum_{M_j^{1/3} < m \le N_j^{1/3}} \!\!\!\!\! \!\! e(a_j\lambda m^3) = \int_{M_j^{1/3}}^{N_j^{1/3}} e(a_j \lambda t^3) d [t] \nonumber \\
& = & \int_{M_j^{1/3}}^{N_j^{1/3}} e(a_j \lambda t^3) dt - \int_{M_j^{1/3}}^{N_j^{1/3}} e(a_j \lambda t^3) d\{t\},
\end{eqnarray}
where $t = [t] + \{t\}$, $[t]$ and $\{t\}$ is the integral and fractional parts of $t$, respectively.

The second term of the right-hand side of \eqref{eq:computeI0} can be bounded by using integration by parts,
\begin{equation} \label{eq:ps1}
\int_{M_j^{1/3}}^{N_j^{1/3}} e(a_j \lambda t^3) d\{t\} = \{t\}e(a_j\lambda t^3)\bigg|_{M_j^{1/3}}^{N_j^{1/3}} - \int_{M_j^{1/3}}^{N_j^{1/3}} \{t\}de(a_j\lambda t^3).
\end{equation}
The first term of the right-hand side of \eqref{eq:ps1} is $\ll 1$, while the second term can be bounded, 
\begin{eqnarray*}
\int_{M_j^{1/3}}^{N_j^{1/3}} \{t\} de(a_j\lambda t^3) & = & 3a_j \lambda \int_{M_j^{1/3}}^{N_j^{1/3}} \{t\} t^2 e(a_j\lambda t^3) dt \\
& \ll & |a_j| |\lambda| \int_{M_j^{1/3}}^{N_j^{1/3}} t^2 dt \\ & \ll & |a_j||\lambda|N_j = N|\lambda|.
\end{eqnarray*}
For the first term of the right-hand side of \eqref{eq:computeI0}, after substitution becomes
\begin{eqnarray*}
\int_{M_j^{1/3}}^{N_j^{1/3}} e(a_j\lambda t^3) dt & = & \frac{1}{3} \int_{M_j}^{N_j} \frac{e(a_j \lambda t)}{t^{2/3}} dt \\
& = & \frac{1}{3} \sum_{M_j < m \le N_j} \frac{e(a_j \lambda m)}{m^{2/3}}.
\end{eqnarray*}
Hence, 
\begin{equation} \label{eq:estimateI0} V_j (\lambda) = \frac{1}{3} \sum_{M_j < m \le N_j} \frac{e(a_j \lambda m)}{m^{2/3}} + O(1 + |\lambda| N). \end{equation}
We also have the elementary bound
\begin{equation} \label{eq:estimateI0main}
\sum_{M_j < m \le N_j} \frac{e(a_j \lambda m)}{m^{2/3}} \ll \min\{N_j^{1/3}, M_j^{-2/3}|a_j \lambda|^{-1}\} \ll |a_j|^{-1/3} \min\{N^{1/3}, {M^{-2/3}} |\lambda |^{-1}\}.
\end{equation}
Therefore, in view of \eqref{eq:estimateI0} and \eqref{eq:estimateI0main},
\begin{equation*}
V_j(\lambda) \ll N_j^{1/3} \min\left\{ 1, \frac{1}{N|\lambda|} \right\} + (1 + |\lambda| N).
\end{equation*}
We can substitute $V_j(\lambda)$ into $H_j$,
\begin{align*}
H_j^2 & = \int_{-1/Q}^{1/Q} |V_j(\lambda)|^2 d\lambda \\
& \ll \int_{-1/N}^{1/N} N_j^{2/3} d\lambda  + \int_{1/N}^{1/Q} N_j^{2/3} \frac{1}{(N|\lambda|)^2} d\lambda + \int_{-1/Q}^{1/Q} (1 + |\lambda|N)^2 d\lambda \\
& \ll |a_j|^{-2/3} N^{-1/3} + |a_j|^{-2/3} N^{-1/3} + \frac{P^3}{N} \\
& \ll |a_j|^{-2/3} N^{-1/3},
\end{align*}
provided that $N_j \ge P^{4.5}$.
\end{proof}

For $\nu = 8,\dots, 1$, we use Cauchy's inequality and Lemma 3.2.6 to obtain
\begin{align} \label{eq:Lnu}
|I_\nu| & \ll  \sum_{\substack{I \subseteq \{1, \dots, 9\} \\ |I| = 2 \\ L = \{1,\dots, 9\} \backslash I}}\left(\prod_{i\in I} K_i\right)\left(\prod_{m \in L} |a_m|^{-1/3} \right) N^{7/3} \nonumber \\
& \ll  \frac{N^2}{|a_1\cdots a_9|^{1/3}L^A}.
 \end{align}

Therefore, it follows that
\begin{equation} |I_1|, \dots, |I_9| \ll \frac{N^2}{|a_1\cdots a_9|^{1/3}L^A}. \label{*}
\end{equation}

In view of (3.20), it remains to compute the main term $I_0$. 
Substituting \eqref{eq:computeI0} into $I_0$ yields
\begin{eqnarray} \label{eq:I0eval}
I_0 & = & \frac{1}{3^9} \sum_{q \le P} \frac{B(n, q)}{\phi^9(q)} \int_{-1/(qQ)}^{1/(qQ)} \prod_{j = 1}^9 \Bigg\{ \sum_{M < |a_j| m \le N} \frac{e(a_j \lambda m)}{m^{2/3}}\Bigg\} e(-n\lambda) d \lambda \nonumber \\
& & + O\Biggr(\sum_{q \le P} \frac{|B(n, q)|}{\phi^9(q)}\int_{-1/(qQ)}^{1/(qQ)} \sum_{\substack{L \subseteq \{1, 2, \dots, 9\} \\ |L| = 8}}  \Biggr\{ \prod_{j \in L} \Bigg|\sum_{M < |a_j|m \le N}\frac{e(a_j\lambda m)}{m^{2/3}} \Bigg| (1 + |\lambda | N) \Biggr\} d\lambda\Biggr) \nonumber. 
\end{eqnarray}

\noindent For the error term, by \eqref{eq:estimateI0} and Lemma \ref{lem:convbdd} with $r_0 = 1$ (because $\chi$ here are principal), we have
\begin{eqnarray} \label{eq:I0bigObdd}
 & & \sum_{q \le P} \frac{|B(n,q)|}{\phi^9(q)} \int_{-1/qQ}^{1/(qQ)}  \Bigg( \Bigg|\sum_{M < |a_j|m \le N}\frac{e(a_j\lambda m)}{m^{2/3}} \Bigg|^8(1 + |\lambda | N) \Bigg) d\lambda \nonumber \\
 && \ll \frac{1}{|a_j|^{8/3}} \sum_{q \le P} \frac{|B(n,q)|}{\phi^9(q)} \left\{ \int_0^{1/(M^{2/3}N^{1/3})} N^{8/3} d \lambda + \int_{1/(M^{2/3}N^{1/3})}^{1/Q} M^{-8/3} \lambda^{-3} d\lambda \right\} \nonumber \\
 && \ll \frac{1}{|a_j|^{4/3}} \left\{ N^{4/3} \lambda \bigg|_0^{1/(M^{2/3}N^{1/3})} + \frac{N}{M^{8/3}} \frac{1}{\lambda^2} \bigg|_{1/(M^{2/3}N^{1/3})}^{1/Q} \right\} \ll \frac{N^{1/3}}{|a_j|^{4/3}}. 
\end{eqnarray}
So by H\"older's inequality,
\begin{eqnarray*}
&& \sum_{q \le P} \frac{|B(n, q)|}{\phi^9(q)}\int_{-1/(qQ)}^{1/(qQ)} \sum_{\substack{L \subseteq \{1, 2, \dots, 9\} \\ |L| = 8}}  \Bigg( \prod_{j \in L} \Bigg|\sum_{M < |a_j|m \le N}\frac{e(a_j\lambda m)}{m^{2/3}} \bigg| (1 + |\lambda | N) \Bigg) d\lambda \\
&& \ll \sum_{q \le P} \frac{|B(n, q)|}{\phi^9(q)}\int_{-1/(qQ)}^{1/(qQ)} \sum_{\substack{L \subseteq \{1, 2, \dots, 9\} \\ |L| = 8}}   \Bigg( \prod_{j \in L} \Bigg|\sum_{M < |a_j|m \le N}\frac{e(a_j\lambda m)}{m^{2/3}} \Bigg| (1 + |\lambda | N)\Bigg)d\lambda \\
&& \ll \sum_{q \le P} \frac{|B(n, q)|}{\phi^9(q)}\sum_{\substack{L \subseteq \{1, 2, \dots, 9\} \\ |L| = 8}} \prod_{j \in L} \Bigg( \int_{-1/(qQ)}^{1/(qQ)}  \Bigg|\sum_{M < |a_j|m \le N}\frac{e(a_j\lambda m)}{m^{2/3}} \Bigg|^8 (1 + |\lambda | N) \Bigg)^{1/8}d\lambda  \\
&& \ll \sum_{\substack{L \subseteq \{1, 2, \dots, 9\} \\ |L| = 8}} \prod_{j \in L} \Bigg(\sum_{q \le P} \frac{|B(n, q)|}{\phi^9(q)} \int_{-1/(qQ)}^{1/(qQ)}  \Bigg|\sum_{M < |a_j|m \le N}\frac{e(a_j\lambda m)}{m^{2/3}} \Bigg|^8 (1 + |\lambda | N)\Bigg)^{1/8} d\lambda. 
\end{eqnarray*}
Using \eqref{eq:I0bigObdd} the last bound becomes
\begin{eqnarray*}
 \ll \frac{N^{1/3}}{|a_1 \cdots a_8|^{1/3}}.
\end{eqnarray*} 
The other error terms in \eqref{eq:I0eval} can be treated similarly and they are $\ll |a_1 \cdots a_9|^{-1/3}N^{2/3}.$ 

We can extend the integral in the main term of \eqref{eq:I0eval} to $[-1/2, 1/2]$; by Lemma \ref{lem:convbdd} and \eqref{eq:estimateI0}, the resulting error is
\begin{eqnarray*}
\ll \frac{1}{|a_1 \cdots a_9|^{1/3}} \int_{1/(PQ)}^{1/2} M^{-6}|\lambda|^{-9} d\lambda \ll \frac{(PQ)^8}{|a_1\cdots a_9|^{1/3} M^6} \ll \frac{N^2}{|a_1\cdots a_9|^{1/3}L^c},
\end{eqnarray*}
where we have used \eqref{eq:PQ}. Thus \eqref{eq:I0eval} becomes
\begin{equation} \label{eq:I0final}
I_0 = \frac{1}{3^9} \SS(n, P)\SI(n) + O\left(\frac{N^2}{|a_1\cdots a_9|^{1/3}L^A}\right).
\end{equation}
Therefore, Lemma \ref{lem:singularbdd} follows from \eqref{eq:split_terms}, \eqref{eq:Lnu}, and \eqref{eq:I0final}. 

This completes the proof of Lemma 3.2.4.
\end{proof}


\subsection{The Singular Series and Singular Integrals}

In this section, we study the singular series and singular integrals.
\begin{lem}
For $j = 1, \dots, 9$, let $\chi_j \!\! \pmod{p^{\alpha_j}}$ be primitive characters and take $\alpha = \max\{\alpha_1, \dots, \alpha_9\}$. For any $t \ge \alpha$ and for the function $B(p^t, \chi_1\chi_0, \dots, \chi_9\chi_0)$ with $\chi_0$ is modulo $p^t$, we have the following:\\ 
\emph{(a)} $B(p^\alpha, \chi_1\chi_0, \dots, \chi_9\chi_0) = F(p^\alpha, \chi_1\chi_0, \dots, \chi_9\chi_0)$, \\ 
\emph{(b)} $B(p^t) = 0$ if $t \ge \theta + \max\{\theta, \alpha\}$, where $\theta = 1$ if $p \ne 3$ and $\theta = 2$ if $p = 3$, \\ 
\emph{(c)} ${\displaystyle \sum_{\nu = \alpha}^\eta \phi(p^\nu)^{-9}B(p^\nu) = \phi(p^\eta)^{-9}F(p^\eta)}$ for any $\eta \ge \alpha.$ 
\end{lem}
\begin{proof}
(a) Without loss of generality, assume $\alpha = \alpha_1 \ge 1$. By comparing \eqref{eq:Bq} and \eqref{eq:Fq} 
it suffices to show that ${\displaystyle \prod_{j=1}^5 C_{\chi_j\chi_0}(a_jh) = 0}$ for each $h$ divisible by $p$. 
However, since $\chi_1$ (mod $p^\alpha$) is primitive, Lemma 2.2.2 (a) asserts that $C_{\chi_1\chi_0}(a_1h) = 0$. This proves (a). \\

(b) This follows directly from (2.4) with $q = p^t$ and Lemma 2.2.2 (b), since by (1.3), $p$ does not divide some $a_j h$. \\

(c) Rewrite the sum in $B(p^\nu)$ as ${\displaystyle \sum_{h=1}^{p^\nu} - \sum_{\substack{h=1 \\ p | h}}^{p^\nu}}$. The first sum is exactly $F(p^\nu)$. When we set $h = ph'$, we can factor $p$ out of each $C_{\chi_j}$ and so  the second sum is precisely $p^9F(p^{\nu-1})$ when $v \ge \max \{ \alpha + 1, 2\}$. So
\[ \phi(p^\nu)^{-9} B(p^\nu) = \phi(p^\nu)^{-9} F(p^\nu) - \phi(p^{\nu-1})^{-9} F(p^{\nu-1}), \]
for $\nu \ge \alpha + 1$. The validity of this relation for $\nu = 1, \alpha = 0$ can be verified directly. By summing both sides for $\nu = \alpha + 1, \dots, \eta$ and using (a) we obtain (c). \end{proof}

\begin{cor} \label{lem:corforA} {\emph{}} \\
\emph{(a)} $A(p^\alpha) = 0$ for primes $p \ne 3, \alpha \ge 2$ and $A(3^\alpha) = 0$ for $\alpha \ge 4$. \\
\emph{(b)} $p^\eta \phi (p^\eta)^{-9}N(p^\eta) = p \phi(p)^{-9} N(p)$ for primes $p \ne 3, \eta \ge 1$. \\
\emph{(c)} $3^\eta \phi(3^\eta)^{-9}N(3^\eta) = 3^3 \phi(3^3)^{-9} N(3^3)$ for $\eta \ge 3$.
\end{cor}
\begin{proof}
To prove (a), take $\chi_1 = \cdots = \chi_9 = \chi_0$ and $\alpha = 0$, as in Lemma 2.2.4. 
Then by \eqref{eq:FqNq} and \eqref{eq:Aq}, we see that $A(p^\nu) = \phi(p^\nu)^{-9}B(p^\nu)$ and $N(p^\nu) = p^{-\nu}F(p^\nu)$. The corollary follows from Lemma 3.4 (b) and (c) since $\theta = 1$ if $p \ne 3$ and $\theta = 3$ if $p = 3$. \end{proof}

\begin{lem} \label{lem:Aineq}
We have $|A(n,p)| < c_1p^{-9/2}$ for all $p \nmid a_1\cdots a_9$, and $p \ne 3$, for some constant \mbox{$c_1$}.
\end{lem}
\begin{proof} From Lemma \ref{lem:chineq} (c), with $k = 1,\dots, p-1$, for each character sum with the principal character,
\begin{eqnarray*}
|C_{\chi_0}(a_jk)| & \le & 3(3,p)(a_jk, p)^{1/2}p^{1/2}  \\
& \ll & p^{1/2}.
\end{eqnarray*}
Then
\begin{eqnarray*}
|A(n,p)| & \le & \phi(p)^{-9}\sum_{k=1}^p \left( \prod_{j=1}^9 |C_{\chi_0}(a_jk)|\right) \\
& \ll & p^{-9} p^{9/2} = p^{-9/2}.
\end{eqnarray*} \end{proof}
\begin{lem}
\emph{(i)} For $x > 0$ and some constant $c_2 > 0$,
\[ \sum_{q > x} |A(n,q)| \ll x^{-1}\log^{c} (x + 2). \]
So the singular series $\SS(n) := \SS(n, \infty)$ is absolutely convergent. \\
\emph{(ii)} We have $\SS(n) \gg (\log\log D)^{-c}$ for some constant $c > 0$. \end{lem}
\begin{proof}
(i) Let $\sigma = (\log(x + 2))^{-1}$. From Lemma \ref{lem:singularbdd} and Corollary \ref{lem:corforA} (a), 
we have
\begin{eqnarray}
\sum_{q > x} |A(n, q)| & \le & \sum_{q = 1}^\infty \left(\frac{q}{x}\right)^{1 - \sigma} |A(n,q)| = x^{-1+\sigma} \sum_{q = 1}^\infty q^{1 - \sigma} |A(n,q)| \nonumber \\
\label{eq:Aineq1} & \ll & x^{-1} \prod_p (1 + p^{1- \sigma} |A(n, p)|), \end{eqnarray}
because $x^\sigma \ll 1$.

From Lemma \ref{lem:Aineq}, 
\begin{eqnarray}
\prod_{p \nmid a_1\cdots a_9} (1+p^{1-\sigma}|A(n,p)|) & \le & x^{-1} \prod_{p\nmid a_1\cdots a_9} \left(1 + \frac{c_1}{p^{4+\sigma}}\right) \nonumber \\
& \le & \prod_p (1-p^{-4-\sigma})^{-c_1} \nonumber \\
\label{eq:Aineq2} & = & \zeta(1 + \sigma)^{c_1} \ll \sigma^{-c_1} = \log^{c_1}(x + 2),
\end{eqnarray}
for some constant $c_1 > 0$. Similarly, for $p |a_1\cdots a_9$, by \eqref{eq:techcond},
\begin{eqnarray}
\prod_{p|a_1\cdots a_9}(1 + p^{1-\sigma}|A(n,p)|) & \le & x^{-1} \prod_{p|a_1\cdots a_9}\left(1 + \frac{c_2}{p^{3+\sigma}}\right) \nonumber\\
\label{eq:Aineq3} & \ll & \log^{c_2}(x+2),
\end{eqnarray}
for some constant $c_2>0$. (i) then follows from \eqref{eq:Aineq1}, \eqref{eq:Aineq2} and \eqref{eq:Aineq3}. \\
(ii) From \eqref{eq:newNq} and \eqref{eq:techcond}, we have $N(p) = p^8 + O(p^7)$. It follows from Lemma \ref{lem:Aineq} that for some large constant $c > c_1$,
\begin{eqnarray*}
\SS(n) & = & \prod_p (1 + A(n,p)) \gg \prod_{\substack{p|a_1\cdots a_9 \\ p > c}} \left(1 - \frac{c}{p}\right) \prod_{\substack{p\nmid a_1\cdots a_9 \\ p > c}} \left( 1 - \frac{c_1}{p^4}\right) \\
& \gg & \prod_{\substack{p|a_1\cdots a_9 \\ p > c}} \left(1 - \frac{c}{p}\right) \gg \prod_{p |a_1\cdots a_9}(1 + p^{-1})^{-(1+ c)}.
\end{eqnarray*}
 Therefore, (ii) follows from the well-known estimate
\[ \prod_{p | x} (1 + p^{-1}) \ll \log \log x. \]
\end{proof}
\begin{lem}
Suppose \eqref{eq:techcond} and either \\
\emph{(i)} $a_j$'s are not all of the same sign and $N \ge C|n|$, for some constant $C$; or \\
\emph{(ii)} all $a_j$'s are positive and $n = N$. \\
Then we have
\begin{equation} \label{eq:SIorder}
\SI(n) = \sum_{\substack{a_1m_1 + \cdots + a_9m_9 = n \\ M < |a_j|m_j^3 \le N}} (m_1 \cdots m_9)^{-2/3} \asymp \frac{N^2}{|a_1 \cdots a_9|^{1/3}}.
\end{equation}
\end{lem} 
\begin{proof} We derive the following inequalities:
\begin{eqnarray*}
\sum_{\substack{a_1m_1 + \cdots + a_9 m_9=n \\ M < |a_j| m_j \le N}} 1 & \le & \sum_{\substack{n - (a_1m_1 + \cdots + a_8m_8) \equiv 0 {\textrm{ (mod $|a_9|$)}} \\ M < |a_j|m_j \le N, j = 1, \dots, 8}} 1 \\
& = & \sum_{\substack{M_j < m_j \le N_j \\ j = 1,\cdots, 7}} \Bigg( \sum_{\substack{m_8 \equiv \overline{a_8}(n - (a_1m_1 + \cdots + a_7m_7)) \textrm{ (mod $|a_9|$)} \\ M_8 < m_8 \le N_8}} 1 \Bigg) \\
& \ll & N_1 N_2 \cdots N_7 \frac{N_8}{|a_9|} \ll \frac{N^8}{|a_1 \cdots a_9|},
\end{eqnarray*}
where $a_8\overline{a_8} \equiv 1$ (mod $|a_9|$).

To establish inequalities in the other direction, we first consider case (ii) in which all $a_j$'s are positive and $n = N$. If $M < a_j m_j \le N/9$ for $j = 1, \dots, 8$, then
\[ M < N/9 = N - 8(N/9) \le N - (a_1m_1 + \cdots + a_8m_8) = a_9 m_9 < N. \]
It follows that
\[ \sum_{\substack{a_1m_1 + \cdots + a_9m_9 = n \\ M < a_j m_j \le N, \; j = 1,\dots, 9}} 1 \ge \sum_{\substack{n - (a_1m_1 + \cdots + a_8m_8) \equiv 0 {\textrm{ (mod $a_9$)}} \\ M < a_j m_j \le N/9, \; j=1,\dots, 8}} 1 \gg \frac{N^8}{|a_1\cdots a_9|}. \]
Case (i) can be treated similarly, and so we have
\[ \sum_{\substack{a_1m_1 + \cdots +a_9m_9 = n \\ M < |a_j| m_j \le N}} 1 \asymp \frac{N^8}{|a_1\cdots a_9|}, \]
 Therefore,
\begin{eqnarray*}
\sum_{\substack{a_1m_1 + \cdots + a_9m_9 = n \\ M < |a_j|m_j \le N}} (m_1 \cdots m_9)^{-2/3} & \asymp & \sum_{\substack{a_1m_1 + \cdots +a_9m_9 = n \\ M < |a_j|m_j \le N}} (N_1 \cdots N_9)^{-2/3} \\
& = & \Bigg(\frac{N^9}{|a_1\cdots a_9|}\Bigg)^{-2/3} \sum_{\substack{a_1m_1 + \cdots+ a_9m_9 = n \\ M < |a_j|m_j \le N}} 1 \\
& \asymp & \Bigg(\frac{|a_1\cdots a_9|^{2/3}}{N^6}\Bigg) \Bigg(\frac{N^8}{|a_1\cdots a_9|}\Bigg) \\
& = & \frac{N^2}{|a_1\cdots a_9|^{1/3}},
\end{eqnarray*}
from which the desired result follows. 
\end{proof}

\subsection{The proof of Lemma 3.2.5} 
In this section, we prove the estimates to $K_j$ and $J_j$ in Lemma \ref{lem:KJlem}. \\

\noindent{\it Proof of Lemma \ref{lem:KJlem}.} Define
\[ \hat{W}_j(\chi, \lambda):=\sum_{M < |a_j|n^3 \le N} \!\!\!\! \Lambda (n) \chi(n) 
e(a_jn^3\lambda) - \delta_\chi \!\!\!\!\! \sum_{M < |a_j|m^3 \le N} e(a_j m^3 \lambda).
\]
Note that
\begin{equation} \label{eq:WjDiff}
W_j(\chi, \lambda) - \hat{W}_j(\chi, \lambda) = -\sum_{m \ge 2} \sum_{M < |a_j|p^{3m} \le N} (\log p ) \chi(p) e(a_jp^{3m} \lambda) \ll N_j^{1/6}. 
\end{equation}
By a dyadic argument: For $1 \le R \le P$, if
\[ g(r) =  {\sum_{\chi \!\!\!\!\! \pmod{r}}}^{\!\!\!\!\!*} \left( \int_{-1/(rQ)}^{1/(rQ)} |W_j(\chi, \lambda)|^{2} d \lambda\right)^{1/2}, \]
then
\begin{align*}
\sum_{r \le P} r^{-1/3+\epsilon} g(r) & = \sum_{0 \le m \le \frac{\log P}{\log 2}} \sum_{\frac{P}{2^{m+1}} \le r \le \frac{P}{2^m}} r^{-1/3+\epsilon}g(r) \\
& \ll \sum_{0 \le m \le \frac{\log P}{\log 2}} \left(\frac{P}{2^m}\right)^{-1/3+\epsilon} \sum_{r \sim \frac{P}{2^m}} g(r) \\
& \ll  \sum_{r \sim \frac{P}{2^m}} g(r) P^{-1/3+\epsilon} \sum_{0 \le m \le \frac{\log P}{\log 2}} 2^{m(1/3-\epsilon)} \\
& \ll  \sum_{r \sim \frac{P}{2^m}} g(r) \sum_{0 \le m \le \frac{\log P}{\log 2}} 1 \\
& \ll \log P  \sum_{r \sim \frac{P}{2^m}} g(r).
\end{align*}
It follows that by the definition of $K_j$ in \eqref{eq:Keq}, we have
\begin{align*}
K_j &\ll  L \max_{R \le P}  \sum_{r \sim R} r^{-1/3+\epsilon} {\sum_{\chi \!\!\!\!\! \pmod{r}}}^{\!\!\!\!\!*} \left( \int_{-1/(rQ)}^{1/(rQ)} |W_j(\chi, \lambda)|^{2} d \lambda\right)^{1/2} \\
& \ll  L \max_{R \le P}  \sum_{r \sim R} r^{-1/3+\epsilon} {\sum_{\chi \!\!\!\!\! \pmod{r}}}^{\!\!\!\!\!*} \left( \int_{-1/(rQ)}^{1/(rQ)} |\hat{W}_j(\chi, \lambda)|^{2} d \lambda\right)^{1/2} + \frac{R^{7/6+\epsilon}L^{c/2} P^{1/2}}{|a_j|^{1/6}N^{1/3}} \\
& \ll  L \max_{R \le P}  \sum_{r \sim R} r^{-1/3+\epsilon} {\sum_{\chi \!\!\!\!\! \pmod{r}}}^{\!\!\!\!\!*} \left( \int_{-1/(rQ)}^{1/(rQ)} |\hat{W}_j(\chi, \lambda)|^{2} d \lambda\right)^{1/2} + |a_j|^{-1/3} N^{-1/6}L^{-A},
\end{align*}
because $N \ge P^{10+\epsilon}|a_j|$ and by \eqref{eq:WjDiff}.

Thus to establish \eqref{eq:Keq} it suffices to show that for any fixed $A > 0$,
\begin{equation} \label{eq:KL2eq2}
\sum_{r \sim R} {\sum_{\chi \!\!\!\!\! \pmod{r}}}^{\!\!\!\!\!*} \left( \int_{-1/(rQ)}^{1/(rQ)} |\hat{W}_j(\chi, \lambda)|^{2} d \lambda\right)^{1/2}  \ll_A |a_j|^{-1/3} R^{1/3-\epsilon}N^{-1/6} L^{-A}
\end{equation}
holds for $R \le P$.

We will need Gallagher's lemma (Lemma 1 in \cite{Gallagher}).
\begin{lem}[Gallagher] \label{lem:Gal} Let 
\begin{equation} \label{eq:Geq1} S(t) = \sum c(\nu) e(\nu t) \end{equation}
be an absolutely convergent exponential sum, and $c(\nu)$ be arbitrary complex numbers. Also, let $\delta = \theta/T$, with $0 < \theta < 1$. Then
\[ \int_{-T}^T |S(t)|^2 dt \ll_\theta \int_{-\infty}^\infty \bigg|\delta^{-1} \sum_{x}^{x + \delta} c(\nu)\bigg|^2 dx. \]
\end{lem}
In view of Gallagher's Lemma, we have
\begin{eqnarray}
\int_{-1/(rQ)}^{1/(rQ)} |\hat{W}_j (\chi, \lambda)|^2d\lambda &\ll & \frac{1}{(RQ)^2} \int_{-\infty}^{\infty} \bigg|\sum_{\substack{t < |a_j|m^3 \le t + rQ \\ M < |a_j|m^3 \le N}} (\Lambda(m)\chi(m) -\delta_\chi)\bigg|^2 dt \nonumber \\ \label{eq:galeq1}
& \ll &  \frac{1}{(RQ)^2} \int_{M-rQ}^{N} \bigg|\sum_{X < m^3 \le Y} (\Lambda(m)\chi(m) -\delta_\chi)\bigg|^2 dt,
\end{eqnarray}
where $X = \max(t, M)/|a_j|$, and $Y =\min(t + rQ, N)/|a_j|$.

Let $D_1, \dots, D_{10}$ be positive numbers such that 
\[ M_j^{1/3} \ll D_1D_2 \cdots D_{10} \ll N_j^{1/3}, \quad {\textrm{and}} \quad D_\nu \le N_j^{1/5}, \quad {\textrm{for $\nu =6,\dots, 10$}}. \]
We also let
\begin{eqnarray*}
a_\nu (n) & := & \left\{\begin{array}{ll} \log n \quad & {\textrm{if $\nu = 1$}}, \\ 1 \quad  & {\textrm{if $1 < \nu \le 5$}}, \\ \mu(n) \quad & {\textrm{if $6 \le \nu \le 10$.}} \end{array}\right.
\end{eqnarray*}
We define the following functions of a complex variable $s$: 
\[ f_j(s,\chi) := \sum_{n \sim D_i} \frac{a_\nu(n)\chi(n)}{n^s}, \qquad F(s,\chi) := \prod_{\nu=1}^{10} f_j (s,\chi).\]
We will state Heath-Brown's identity (see \S4 of Choi \cite{Choi1}) for $k = 5$, which says that
\[ \frac{\zeta'}{\zeta}(s) = \sum_{\nu = 1}^5 \binom{5}{\nu}(-1)^{\nu-1}\zeta'(s)\zeta^{\nu-1}(s) + \frac{\zeta'}{\zeta}(s)(1-\zeta(s)G(s))^5, \]
where $\zeta(s)$ is the Riemann zeta-function, and 
\[ G(s) = \sum_{m \le N_j^{1/15}}\mu(m) m^{-s}. \]
The reason why we choose $k = 5$ is that the identity with $k \le 4$ will give weaker results, and when $k \ge 6$ it produces the same estimate as the case $k = 5$. Equating coefficients of the Dirichlet series on both sides provides an identity for $-\Lambda(m)$. Also, for $m \le N_j^{1/3}$ the coefficient of $m^{-s}$ in $-(\zeta'/\zeta)(s)(1-\zeta(s)G(s))^5$ is zero. Thus,
\[ \Lambda(m) = \sum_{\nu=1}^5 \binom{5}{\nu}(-1)^{\nu-1} \sum_{\substack{m_1\cdots m_{2\nu} =m\\ m_{\nu+1}\cdots m_{2\nu} \le N_j^{1/3}}} (\log m_1) \mu(m_{\nu+1}) \cdots \mu(m_{2\nu}). \]
Applying this identity to the inner sum in \eqref{eq:galeq1}, 
\begin{eqnarray} \label{eq:Lamchi}
 \sum_{X < m^3 \le Y} (\Lambda(m)\chi(m) - \delta_\chi), 
\end{eqnarray}
by the dyadic argument again, we find that \eqref{eq:Lamchi} is a linear combination of $O(L^{10})$ terms, each of which is of the form 
\[ \sigma(t; {\bf D}) = \underset{\substack{ \ \ \  m_1 \sim D_1 \  m_{10} \sim D_{10} \\ M_j^{1/3} < m_1 \cdots m_{10} \le N_j^{1/3}}}{\sum \, \cdots \, \sum} a_1(m_1)\chi(m_1) \cdots a_{10}(m_{10})\chi(m_{10}), \]
where {\bf D} denotes the vector $(D_1, \dots, D_{10})$. We need Perron's summation formula (see Lemma 3.12 in \cite{Titchmarsh1} for example).
\begin{thm2}[Perron's formula] Let ${\displaystyle f(s) = \sum_{n=1}^\infty \frac{a_n}{n^s}}$ be a Dirichlet series that is absolutely convergent for {\emph{Re}}$(s) >1$. For $x$ not an integer and $\sigma > 1$, we have
\[ \sum_{n \le x} a_n = \frac{1}{2\pi i} \int_{\sigma-iT}^{\sigma+iT} f(s) \frac{x^s}{s} ds + O\bigg(\sum_{n=1}^\infty \left(\frac{x}{n}\right)^\sigma |a_n| \min\left(1, \frac{1}{T|\log(x/n)|}\right)\bigg). \]
\end{thm2}
By using Perron's formula and then shifting the contour to the left, the above $\sigma(u; {\bf D})$ is
\begin{eqnarray*}
& = & \frac{1}{2\pi i}\int_{1 + 1/L - iT}^{1 + 1/L + iT} F(s, \chi) \frac{Y^{s} -X^{s}}{s} ds + O\left(\frac{N_j^{1/3}L^2}{T}\right) \\
& = & \frac{1}{2\pi i}\left\{\int_{1 + 1/L -iT}^{1/2 - iT} + \int_{1/2-iT}^{1/2+iT} + \int_{1/2+iT}^{1+1/L+iT}\right\} + O\left(\frac{N_j^{1/3}L^2}{T}\right),
\end{eqnarray*}
where $T$ is a parameter satisfying $2 \le T \le N_j^{1/3}$. The integral on the two horizontal segments above can be easily estimated as
\[ \ll \max_{1/2 \le \sigma \le 1 + 1/L} |F(\sigma \pm iT, \chi)|\frac{u^\sigma}{T} \ll \max_{1/2 \le \sigma \le 1 + 1/L} N_j^{(1-\sigma)/3}L\frac{Y^\sigma}{T} \ll \frac{N_j^{1/3}L}{T}, \]
on using the trivial estimate
\[ F(\sigma \pm iT, \chi) \ll \prod_{\nu =1}^{10} |f_\nu(\sigma \pm iT, \chi)| \ll (D_1^{1 -\sigma}L) \prod_{\nu=2}^{10}D_\nu^{1-\sigma} \ll N_j^{(1-\sigma)/3}L. \]
Thus,
\begin{eqnarray} \label{eq:sigmaTD}
\sigma(t; {\bf D}) = \frac{1}{2\pi i}\int_{-T}^T F\left(\frac{1}{2}+iu, \chi\right)\frac{Y^{\frac{1}{3}\left(\frac{1}{2}+iu\right)} - X^{\frac{1}{3}\left(\frac{1}{2}+iu\right)}}{\frac{1}{2}+iu}du + O\left(\frac{N_j^{1/3}L^2}{T}\right). \qquad
\end{eqnarray}
Note that for any $0 < \beta < 1$, 
\begin{eqnarray} \label{eq:YXbdd}
Y^\beta - X^\beta \ll \frac{(t + rQ)^\beta - t^\beta}{|a_j|^\beta} = \frac{t^\beta\{(1+rQ/t)^\beta - 1\}}{|a_j|^\beta} \ll \frac{rQ}{|a_j|^\beta M^{1-\beta}},
\end{eqnarray}
where in the last step we used $M - rQ \le t \le N$ and $rQ \le 2RQ \le 2 PQ \ll ML^{-c}$, and the fact that
\begin{eqnarray*}
t^\beta\{(1+rQ/t)^\beta - 1\} & \ll & t^\beta\{(1+rQ/t)-1\} = t^\beta(rQ/t) = \frac{rQ}{t^{1-\beta}} \\
& \ll & \frac{rQ}{(M-rQ)^{1-\beta}} = \frac{rQ}{M^{1-\beta}}\left(1 + \frac{rQ}{M} + O\left(\frac{rQ^2}{M^2}\right)\right) \ll \frac{rQ}{M^{1-\beta}}.
\end{eqnarray*}
For $\chi = \chi_0$ (mod 1), \eqref{eq:Lamchi} becomes
\[ Y^{1/3} - X^{1/3} \ll |a_j|^{-1/3}M^{-2/3}Q \]
by \eqref{eq:YXbdd} with $r = 1$. This contributes to \eqref{eq:galeq1} acceptably. 

For $\chi \ne \chi_0$ (mod 1), we have $\delta_\chi = 0$ in \eqref{eq:Lamchi}. Then one can see that
\begin{eqnarray*}
\frac{Y^{\frac{1}{3}\left(\frac{1}{2}+iu\right)} - X^{\frac{1}{3}\left(\frac{1}{2}+iu\right)}}{\frac{1}{2}+iu} = \frac{1}{3}\int_X^Y t^{-5/6+iu/3}dt = \frac{1}{3}\int_X^Y t^{-5/6}e\left(\frac{u}{6\pi} \log t\right) dt \ll Y^{1/6} - X^{1/6}.
\end{eqnarray*}
The integral can be easily estimated by \eqref{eq:YXbdd} as $\ll Y^{1/6} - X^{1/6} \ll |a_j|^{-1/6}M^{-5/6}RQ$. On the other hand, one has trivially
\begin{eqnarray*}
\frac{Y^{\frac{1}{3}\left(\frac{1}{2}+iu\right)} - X^{\frac{1}{3}\left(\frac{1}{2}+iu\right)}}{\frac{1}{2}+iu} \ll \frac{Y^{1/6}}{|u|} \ll \frac{N_j^{1/6}}{|u|}.
\end{eqnarray*}
Together with the two upper bounds yields
\begin{eqnarray*}
\frac{Y^{\frac{1}{3}\left(\frac{1}{2}+iu\right)} - X^{\frac{1}{3}\left(\frac{1}{2}+iu\right)}}{\frac{1}{2}+iu} & \ll & \min\left(\frac{RQ}{M^{5/6}|a_j|^{1/6}}, \frac{N_j^{1/6}}{|u|}\right) \ll \frac{1}{|a_j|^{1/6}}\min\left(\frac{RQ}{N^{5/6}}, \frac{N^{1/6}}{|u|}\right).
\end{eqnarray*}
Taking $T = N_j^{1/3}$ and $T_0 = N/QR$, we see that
\begin{align}
\sigma(t; {\bf D})  \ll &  \frac{RQ}{|a_j|^{1/6}N^{5/6}}\int_{|u|\le T_0} \bigg|F\left(\frac{1}{2}+iu, \chi\right)\bigg|du \nonumber \\ \label{eq:Ksigma}
& + \frac{N^{1/6}}{|a_j|^{1/6}}\int_{T_0 < |u| \le T} \bigg|F\left(\frac{1}{2}+iu, \chi\right)\bigg|\frac{dt}{|u|} +L^{12}.
\end{align}

Consequently, from \eqref{eq:Lamchi},
\begin{align*}
 \sum_{X < m^3 \le Y} (\Lambda(m)\chi(m) - \delta_\chi) \ll & \frac{RQL^{10}}{|a_j|^{1/6}N^{5/6}} \max_{\text{\bf D}} \int_{|u| \le T_0} \left|F\left(\frac{1}{2}+iu, \chi\right)\right| du \\ &+ \max_{\text{\bf D}} \frac{N^{1/6}L^{10}}{|a_j|^{1/6}}\int_{T_0 < |u| \le T} \bigg|F\left(\frac{1}{2}+iu, \chi\right)\bigg|\frac{dt}{|u|} + L^{12}.
\end{align*}
Hence from \eqref{eq:galeq1},
\begin{eqnarray*}
\int_{-1/(rQ)}^{1/(rQ)} |\hat{W}_j(\chi, \lambda)|^2d\lambda & \ll & \max_{\bf D} \bigg(\int_{|u|\le T_0} \left|F\left(\frac{1}{2} + iu, \chi\right)\right|du\bigg)^2  \\
&& + \max_{\bf D} \bigg(\int_{T_0 < |u| \le T} \left|F\left(\frac{1}{2}+iu, \chi\right)\right|\frac{du}{|u|}\bigg)^2 \\
& \ll & \frac{L^{20}}{|a_j|^{1/3}N^{2/3}}\max_{\bf D} \bigg(\int_{|u|\le T_0} \left|F\left(\frac{1}{2} + iu, \chi\right)\right|\bigg)^2  \\
&& + \frac{L^{20}N^{4/3}}{|a_j|^{1/3}(RQ)^2}\max_{\bf D} \bigg(\int_{T_0 < |u| \le T} \left|F\left(\frac{1}{2}+iu, \chi\right)\right|\frac{du}{|u|}\bigg)^2 \\
&& + \frac{L^{c}N}{(RQ)^2}.
\end{eqnarray*}
Taking square-root, and summing over primitive characters $\chi \!\!\pmod{r}$ and then summing over $r\sim R$, 
\begin{eqnarray*}
&& \sum_{r \sim R} {\sum_{\chi \!\!\!\!\! \pmod{r}}}^{\!\!\!\!\!*} \bigg(\int_{-1/(rQ)}^{1/(rQ)} |\hat{W}(\chi, \lambda)|^2d\lambda \bigg)^{1/2} \\
& &\ll  \frac{L^{10}}{|a_j|^{1/6}N^{1/3}} \sum_{r \sim R} {\sum_{\chi \!\!\!\! \pmod{r}}}^{\!\!\!\!\!*}  \max_{\bf D} \bigg(\int_{|u|\le T_0} \left|F\left(\frac{1}{2} + iu, \chi\right)\right|\bigg)  \\
&& + \frac{L^{10}N^{2/3}}{|a_j|^{1/6}(RQ)} \sum_{r \sim R} {\sum_{\chi \!\!\!\!\! \pmod{r}}}^{\!\!\!\!\!*}  \max_{\bf D} \bigg(\int_{T_0 < |u| \le T} \left|F\left(\frac{1}{2}+iu, \chi\right)\right|\frac{du}{|u|}\bigg) \\
&& + \frac{L^{c/2}N^{1/2}R}{Q}.
\end{eqnarray*}
Thus to prove \eqref{eq:KL2eq2} it suffices to show that the estimate
\begin{equation} \label{eq:Dpolyeq1}
\sum_{r \sim R} {\sum_{\chi \!\!\!\!\! \pmod{r}}}^{\!\!\!\!\!*} \int_{T_1}^{2T_1} \bigg|F\left(\frac{1}{2}+iu, \chi\right)\bigg| du \ll R^{1/3-\epsilon} N_j^{1/6}
\end{equation}
holds for $R \le P$ and $0 < T_1 \le T_0$, and the estimate
\begin{equation} \label{eq:Dpolyeq2}
\sum_{r \sim R} {\sum_{\chi \!\!\!\!\! \pmod{r}}}^{\!\!\!\!\!*} \int_{T_2}^{2T_2} \bigg|F\left(\frac{1}{2} + iu, \chi\right)\bigg|du \ll R^{1/3-\epsilon} \left(\frac{RQ}{N^{5/6}|a_j|^{1/6}}\right)T_2L^c
\end{equation}
holds for $R \le P$ and $T_0 < T_2 \le T$. 

To prove \eqref{eq:Dpolyeq1} and \eqref{eq:Dpolyeq2} we need the following two lemmas.
\begin{lem} \label{lem:Pan1} For any $P \ge 1$, $T\ge 1$, and $k = 0, 1$,
\[
\sum_{r \le P}  {\sum_{\chi \!\!\!\!\! \pmod{r}}}^{\!\!\!\!\!*} \int_{-T}^T \bigg| L^{(k)}\left(\frac{1}{2}+it, \chi\right)\bigg|^4 dt \ll P^2 T(\log PT)^{4(k+1)}. 
\]
\end{lem}
\begin{lem} \label{lem:Pan2} For any $P \ge 1$, $T \ge 1$ and any complex numbers $a_n$, 
\[
\sum_{r \le P}  {\sum_{\chi \!\!\!\!\! \pmod{r}}}^{\!\!\!\!\!*} \int_{-T}^T \bigg|\sum_{n = M}^{M+N} a_n \chi(n)n^{-it}\bigg|^2 dt \ll \sum_{n = M}^{M+N}(P^2T+n)|a_n|^2. 
\]
\end{lem}
The proofs for Lemma \ref{lem:Pan1} and Lemma \ref{lem:Pan2} can be found in \cite{Pans}, Chapters 2 and 3, respectively. 
\begin{prop} \label{prop:DpolyProp1}
If there exist natural numbers $D_k, D_l$, with $1 \le k,l \le 5$, such that their product $D_{k}D_{l} \ge P^{4/3}$, then the estimate \eqref{eq:Dpolyeq1} holds. 
\end{prop}

\begin{proof} Without loss of generality, suppose that $l =1$, $D_1 = \log n$ and $k = 2$, $D_2 = 1$. Arguing exactly as in the proof of Proposition 1 in Zhan \cite{TZhan1}, we find for $f_1$,
\begin{eqnarray*}
f_1\left(\frac{1}{2}+it, \chi\right) \ll L\bigg(\int_{-N_j^{1/3}}^{N_j^{1/3}}\left|L'\left(\frac{1}{2}+it+iu, \chi\right)\right|^4\frac{du}{1+|u|}\bigg)^{1/4} + L,
\end{eqnarray*}
and so
\begin{eqnarray*}
 && \sum_{r \sim R}  {\sum_{\chi \!\!\!\!\! \pmod{r}}}^{\!\!\!\!\!*} \int_{T_1}^{2T_1}\left|f_1\left(\frac{1}{2}+it, \chi\right)\right|^4 dt \\
 && \ll L^4 \int_{-N_j^{1/3}}^{N_j^{1/3}} \frac{du}{1+|u|} \sum_{r \sim R}  {\sum_{\chi \!\!\!\!\! \pmod{r}}}^{\!\!\!\!\!*} \int_u^{T_1 + u} \left|L'\left(\frac{1}{2}+it, \chi\right)\right|^4 dt + T_1R^2 L^4 \\
 && \ll L^5 \max_{|P| \le N_j^{1/3}} \int_{P/2}^P \frac{du}{1+|u|} \sum_{r \sim R}  {\sum_{\chi \!\!\!\!\! \pmod{r}}}^{\!\!\!\!\!*} \int_u^{T_1 + u} \left|L'\left(\frac{1}{2}+it, \chi\right)\right|^4 dt + T_1 R^2 L^4.
 \end{eqnarray*}
Interchanging order of summation and using Lemma \ref{lem:Pan2}, 
\begin{eqnarray*}
&& \ll L^4 \max_{|P| \le N_j^{1/3}} P^{-1} \int_{T_1}^{2T_1} \sum_{r \sim R}  {\sum_{\chi \!\!\!\!\! \pmod{r}}}^{\!\!\!\!\!*} \int_{(P/2)+t}^{P+t} \left|L'\left(\frac{1}{2}+iu, \chi\right)\right|^4 du\; dt + T_1 R^2 L^4  \\
&& \ll R^2 T_0 L^c,
\end{eqnarray*}
since $0 \le T_1 \le T_0$. The inequality holds for $f_2$ as well, with an extra power of $\log N$. 

Using \ref{lem:Pan2} and H\"older's inequality, we obtain
\begin{eqnarray*}
 && \sum_{r \sim R}  {\sum_{\chi \!\!\!\!\! \pmod{r}}}^{\!\!\!\!\!*} \int_{T_1}^{2T_1} \left|F\left(\frac{1}{2}+it, \chi\right)\right|dt \\
 && \ll \bigg(\sum_{r \sim R}  {\sum_{\chi \!\!\!\!\! \pmod{r}}}^{\!\!\!\!\!*} \int_{T_1}^{2T_1} \left|f_1\left(\frac{1}{2}+it, \chi\right)\right|^4 dt \bigg)^{1/4} \\
  && \  \,  \times \bigg(\sum_{r \sim R}  {\sum_{\chi \!\!\!\!\! \pmod{r}}}^{\!\!\!\!\!*} \int_{T_1}^{2T_1} \left|f_2\left(\frac{1}{2}+it, \chi\right)\right|^4 dt \bigg)^{1/4} \\
    && \ \, \times \bigg(\sum_{r \sim R}  {\sum_{\chi \!\!\!\!\! \pmod{r}}}^{\!\!\!\!\!*} \int_{T_1}^{2T_1} \left|\prod_{k=3}^{10} f_k\left(\frac{1}{2}+it, \chi\right)\right|^2 dt \bigg)^{1/2} \\
    && \ll (R^2T_0)^{1/2} \bigg(R^2T_0 + \frac{N_j^{1/3}}{D_1D_2}\bigg)^{1/2}L^c \\
    && \ll N_j^{1/6}R^{1/3-\epsilon}L^{-A},
\end{eqnarray*}
by the definition of $T_0$ and the condition of the proposition. 

\begin{prop} \label{prop:DpolyProp1} Let $V = \{1, 2, \dots, 10\}$. If $V$ can be divided into two disjoint subsets $V_1$ and $V_2$ such that
\[ \max\bigg\{\prod_{i \in V_1} D_i, \prod_{i \in V_2} D_i \bigg\} \ll N_j^{1/3}P^{-4/3-\epsilon}, \]
then the estimate \eqref{eq:Dpolyeq1} holds. 
\end{prop}
\begin{proof} Denote 
\[ S_k := \prod_{i \in V_k} D_i, \quad k = 1, 2. \]
Also, for $k = 1, 2$, let
\begin{eqnarray*}
 F_k (s, \chi) & = & \prod_{i \in V_k} f_i (s, \chi) \\
 & = & \sum_{n \ll S_k } b_k (n) \chi(n) n^{-s},
\end{eqnarray*}
where $b_k(n)$ is a convolution of the coefficients $a_k(n)$, with the property that {\mbox{$b_k (n) \ll d^c(n)$,}} for some constant $c > 0$. Applying Lemma \ref{lem:Pan2} and the fact that
\[ \sum_{n \le x} d^k(n) \ll x(\log x)^{k(c)}, \]
we have
\begin{eqnarray} \label{eq:KforJ}
\sum_{r \sim R}  {\sum_{\chi \!\!\!\!\! \pmod{r}}}^{\!\!\!\!\!*} \int_{T_0}^{2T_0} \left|F\left(\frac{1}{2}+it, \chi\right)\right|dt & \ll & \prod_{k=1}^2 \sum_{r \sim R}  {\sum_{\chi \!\!\!\!\! \pmod{r}}}^{\!\!\!\!\!*} \int_{T_0}^{2T_0} \left|F_k\left(\frac{1}{2}+it, \chi\right)\right|dt \qquad \\
& \ll & (R^2T_0 +S_1)^{1/2}(R^2T_0 + S_2)^{1/2}.  \nonumber
\end{eqnarray}
If $S_1, S_2 \le N_j^{1/3} P^{-4/3-\epsilon}$, then the above becomes
\begin{eqnarray} \label{eq:beforeP1}
& \ll & R^2T_0 + RT_0^{1/2} N_j^{1/6}P^{-2/3-\epsilon} + N_j^{1/6}L^C.
\end{eqnarray}
 From the proof of the proposition an estimate $P \ll N_j^{1/10 - \epsilon}$ would suffice. \end{proof}

  Now we can finish proving \eqref{eq:Dpolyeq1}. In view of Proposition \ref{prop:DpolyProp1} we may assume
\[ N_j^{1/15} \le D_kD_l \le P^{4/3+\epsilon} \le  N_j^{2/15}, \quad 1 \le k,l \le 5, \quad {\textrm{with $k \ne l$.}} \]
Therefore, by the pigeon-hole principle, there exists at most one $D_k$, with $1 \le k \le 10$ such that $D_k \ge N_j^{1/15}$. If exists, denote it by $D_{k_0}$, otherwise, take $D_{k_0} = 1$. Reorder the remaining $D_k$ as follows:
\[ D_{k_1} \ge D_{k_2} \ge \cdots \ge D_{k_B}, \quad {\text{where }} B = 9 {\text{ or }} 10. \]
Find an integer $1 \le l \le B-1$ such that
\[ \prod_{h=0}^{l-1} D_{k_h} \le N_j^{2/15}, \qquad {\textrm{but}} \qquad \prod_{h=0}^{l} D_{k_h} \ge N_j^{2/15}. \]

Denote 
\[ S_1 := \prod_{h=0}^l D_{k_h}, \qquad {\textrm{and}} \qquad S_2 := \prod_{h=0}^B D_{k_h}. \]
We therefore have
\[ S_1 \ll N_j^{2/15} D_{k_l} \le N_j^{1/5}, \qquad {\textrm{and}} \qquad S_2 \ll N_j^{1/15}S_1^{-1} \ll N_j^{1/5}. \]
The two sets $S_1$ and $S_2$ satisfy the conditions of Proposition \label{prop:DpolyProp2}. Hence \eqref{eq:Dpolyeq1} is proved. \end{proof}

We can now prove (3.26).
We have
\[ J_j \ll L \max_{R \le P} J_j(R), \]
where $J_j(R)$ is defined similarly to $J_j$ except that the sum is over $r \sim R$. The estimation of $J_j(R)$ falls naturally into two cases depending on $R$ is small or large. For $R > L^C$, where $C$ is some positive constant, one can use the machinery that was already developed for $K_j$ in Lemma \ref{lem:KJlem}. We will prove this in Lemma \ref{lem:JLarge}. While for $R \le L^C$, one uses the classical zero-density estimate and zero-free region of the Dirichlet $L$-functions, as we will show in Lemma \ref{lem:JSmall}. 

We first establish the following result for large $R$. 
\begin{lem} \label{lem:JLarge} There exists a constant $c = c(A) > 0$ such that when $L^C < R \le P$, 
\[ J_j(R) \ll N_j^{1/3} L^{-A}, \]
for all $A > 0$, where the implied constant depends at most on $A$. 
\end{lem}
\begin{proof} To prove this result, it suffices to show that
\begin{equation} \label{eq:Jsuffice}
\sum_{r \sim R}  {\sum_{\chi \!\!\!\!\! \pmod{r}}}^{\!\!\!\!\!*} \max_{|\lambda| \le 1(rQ)} |W_j(\chi, \lambda)| \ll R^{1/3-\epsilon} N_j^{1/3}L^{-A}, 
\end{equation}
holds for $L^c < R \le P$ and arbitrary $A$. From \eqref{eq:WjDiff}, it is more convenient to use $|\hat{W}(\chi, \lambda)|$ in place of $|W(\chi, \lambda)|$  with difference $O(N_j^{1/6})$. 

Following the proof of \eqref{eq:Keq}, recall \eqref{eq:sigmaTD},
\begin{eqnarray*}
\sigma(t; {\bf D}) = \frac{1}{2\pi i}\int_{-T}^T F\left(\frac{1}{2}+iu, \chi\right)\frac{t^{\frac{1}{3}\left(\frac{1}{2}+iu\right)} - X^{\frac{1}{3}\left(\frac{1}{2}+iu\right)}}{\frac{1}{2}+iu}du + O\left(\frac{N_j^{1/3}L^2}{T}\right).
\end{eqnarray*}
Using Riemann-Stieljes integration, integration by parts on $\hat{W}_j(\chi, \lambda)$, and then re-arranging the terms, $\hat{W}_j(\chi, \lambda)$ will be a linear combination of $O(L^{10})$ terms, each of which is of the form
\begin{align*}
\hat{W}_j(\chi, \lambda) & = \sum_{M < |a_j|m^3 \le N} \Lambda(m)\chi(m) e(a_jm^3 \lambda) = \int_{M_j^{1/3}}^{N_j^{1/3}} e(a_ju^3\lambda)d\Bigg\{\sum_{M_j^{1/3} < m\le u} \Lambda(m) \chi(m) \Bigg\} \\
& = \frac{1}{2\pi} \int_{-T}^T F\left(\frac{1}{2} + it, \chi\right) \int_{M_j^{1/3}}^{N_j^{1/3}} u^{-1/2+it} e(a_ju^3\lambda) du \; dt + O\left(\frac{N_j^{1/3}L^2}{T}(1+|\lambda|N)\right).
\end{align*}

By taking $T = N_j^{1/3}$ and changing variables in the inner integral, we deduce from the above that
\begin{align} \label{eq:Jjinnerint}
|\hat{W}_j(\chi, \lambda)| & \ll L^{10} \max_{\text{\bf D}} \bigg| \int_{-T}^T F\left(\frac{1}{2} + it, \chi\right) \int_{M_j}^{N_j} v^{-5/6}e\left(t\frac{\log v}{6\pi} + a_j v \lambda\right) dv\, dt \bigg| + PL^c,
\end{align}
where the maximum is taken over all {\bf D} = $(D_1, \dots, D_{10})$. Since
\begin{align*}
\frac{d}{dv} \left(\frac{t}{6\pi} + a_j \lambda v \right) = \frac{t}{6\pi v} + a_j \lambda, \quad \frac{d^2}{dv^2} \left(\frac{t}{6\pi} \log v + a_j \lambda v\right) = -\frac{t}{6\pi v^2},
\end{align*}
by Lemma 4.4 and Lemma 4.3 in \cite{Titchmarsh1}, the inner integral in \eqref{eq:Jjinnerint} can be estimated as
\begin{align} \label{eq:Jjinnerint2}
\ll M_j^{-5/6} \min \bigg \{ \frac{N_j}{(|t| + 1)^{1/2}}, \frac{N_j}{\displaystyle \min_{M_j < v \le N_j} |t + 6\pi a_j \lambda v|}\bigg\} \ll \left\{ \begin{array}{ll} N_j^{1/6} (|t|+1)^{-1/2} & {\textrm{ if $|t| \le T_0$}},  \\ N_j^{1/6}|t|^{-1} & {\textrm{ if $T_0 < |t| \le T$,}} \end{array}\right.
\end{align}
where $T_0 = 8\pi N/(RQ)$. Here, the choice of $T_0$ is to ensure that $|t + 4\pi a_j \lambda u| > |t|/2$ whenever $|t| > T_0$. In fact,
\[ |t + 4\pi a_j \lambda v| \ge |t| - 4\pi |a_j v|/(rQ) > |t|/2 + T_0/2 - 4\pi N/(RQ) = |t|/2. \]
It therefore follows from \eqref{eq:Jjinnerint} and \eqref{eq:Jjinnerint2} that the lemma (more precisely,  \eqref{eq:Jsuffice}) is a consequence of the following two estimates: For $0 < T_1 \le T_0$, we have
\begin{equation} \label{eq:newJint1}
\sum_{r \sim R}  {\sum_{\chi \!\!\!\!\! \pmod{r}}}^{\!\!\!\!\!*} \int_{T_1}^{2T_1} \left|F\left(\frac{1}{2}+it, \chi\right)\right| dt \ll R^{1/3-\epsilon} N_j^{1/6} (T_2+1)^{1/2}L^{-A},
\end{equation}
while for $T_0 < T_2 \le T$, we have
\begin{equation} \label{eq:newJint2}
\sum_{r \sim R}  {\sum_{\chi \!\!\!\!\! \pmod{r}}}^{\!\!\!\!\!*} \int_{T_2}^{2T_2} \left|F\left(\frac{1}{2}+it, \chi\right)\right| dt \ll R^{1/3-\epsilon} N_j^{1/6} T_2 L^{-A}.
\end{equation}
Both \eqref{eq:newJint1} and \eqref{eq:newJint2} are be deduced from the left-hand side of \eqref{eq:beforeP1}. For example, by taking $T_0 =T_1$ we can see from \eqref{eq:KforJ} that
\[ \ll R^2T_1 + RT_1^{1/2} N_j^{1/6}P^{-2/3-\epsilon} + N_j^{1/6}L^c \ll R^{1/3-\epsilon} N_j^{1/6} (T_1+1)^{1/2}L^{-A}, \]
provided $L^c < R \le P = (N/D)^{1/10-\epsilon}$  with a sufficiently large $c$. Here, $L^c < R$ guarantees that the $L^c$ in \eqref{eq:beforeP1} is dominated by the quantity on the right-hand side. This establishes \eqref{eq:newJint1}. Similarly, we can prove \eqref{eq:newJint2} by taking $T_0 = T_2$.  Lemma \ref{lem:JLarge} now follows. \end{proof}

\begin{lem}\label{lem:JSmall} Let $c > 0$ be arbitrary. For $R \le L^c$ and for any fixed $A > 0$, we have
\[ J_j(R) \ll N_j^{1/3}L^{-A}, \]
where the implied constant depends at most on $c$. 
\end{lem}
\begin{proof} We use the explicit formula (see pg. 109, \S17 and pg. 120, \S19 in \cite{Davenport1})
\begin{equation} \label{eq:explicitformula}
\sum_{m \le u} \Lambda(m) \chi(m) = \delta_\chi u - \sum_{|\gamma| \le T} \frac{u^\rho}{\rho} + O\left(\left(\frac{u}{T}+1\right)\log^2(quT)\right),
\end{equation}
where $\beta + i\gamma$ is a non-trivial zero of the function $L(s, \chi)$, and $2 \le T \le u$ is a parameter. Taking $T= N_j^{1/6}$ in \eqref{eq:explicitformula} and inserting it into $\hat{W}_j(\chi, \lambda)$, by $M_j^{1/3} < u \le N_j^{1/3}$, $M_j = CN_j$, and \eqref{eq:Wvlog} we have
\begin{align*}
\hat{W}_j(\chi, \lambda) & = \int_{M_j^{1/3}}^{N_j^{1/3}} e(a_j u^3 \lambda) d\left\{ \sum_{n \le u} (\Lambda(m)\chi(m) - \delta_\chi)\right\} \\
& = -\int_{M_j^{1/3}}^{N_j^{1/3}} e(a_ju^3\lambda) \sum_{|\gamma| \le T} u^{\rho-1} du + O\left(\frac{N_j^{1/3}}{T}(1+|\lambda|N)L^2\right) \\
& = N_j^{1/3}\sum_{|\gamma| \le T} N_j^{(\beta-1)/3} + O\left(\frac{N_j^{1/3}}{T}PL^c\right) \\
& = N_j^{1/3}\sum_{|\gamma| \le N_j^{1/6}} N_j^{(\beta-1)/3} + O(N_j^{1/6}PL^c).
\end{align*}
Now we need Satz VIII.6.2 in Prachar \cite{Prachar}, which states that ${\displaystyle \prod_{\chi \!\! \pmod{q}} L(s, \chi)}$ is zero-free in the region $\sigma \ge 1 - \eta(T)$, $|t| \le T$ except for the possible Siegel zero, where $\eta(T) = c_3 \log^{-4/5}T$. But by Siegel's theorem (see for example \cite{Davenport1}, \S21) the Siegel zero does not exist in the present situation, because $r \le L^C$. We also need the zero-density estimate (see e.g. Huxley \cite{Huxley}):
\[ N^*(\alpha, q, T) \ll (qT)^{12(1-\alpha)/5}\log^c(qT), \]
where $N^*(\alpha, q, T)$ denotes the number of zeros of ${{\displaystyle \prod_{\chi \!\! \pmod{q}}}^{\!\!\!\!\! *} \, L(s, \chi)}$ in the region Re$(s) \ge \alpha$, $|$Im$(s)| \le T$. Thus,
\begin{align*}
\sum_{|\gamma| \le N_j^{1/6}} N_j^{(\beta-1)/3} & \ll L^c \int_0^{1-\eta(N_j^{1/6})} (N_j^{1/6})^{12(1-\alpha)/5}N_j^{(\alpha-1)/3}d\alpha \ll L^c N_j^{-\eta(N_j^{1/6})/10} \\ 
& \ll \exp(-c_4L^{1/5}).
\end{align*}
 Consequently,
 \begin{equation} \label{eq:finalJj2}
 \sum_{r \sim R}  {\sum_{\chi \!\!\!\!\! \pmod{r}}}^{\!\!\!\!\!*} \max_{|\lambda| \le 1/(rQ)} |\hat{W}_j(\chi, \lambda)| \ll N_j^{1/3}L^{-A},
 \end{equation}
 where $R \le L^C$, and $A > 0$ arbitrary. Lemma \ref{lem:JSmall} now follows from \eqref{eq:finalJj2} and \eqref{eq:WjDiff}. \end{proof}

\section{Treatment of the Minor Arcs}
This minor arc can be treated by standard methods such as Hua's and Weyl's inequality, which the detailed treatments can be found in Vaughan \cite{Vaughan}. 

We first derive
\begin{eqnarray} \label{eq:minor0}
\int_\MINA |S_1(\alpha)\cdots S_9(\alpha)|d\alpha & \le & \bigg(\sup_\MINA |S_9(\alpha)|\bigg) \int_0^1 |S_1(\alpha)\cdots S_8(\alpha)|d\alpha.
\end{eqnarray}
We can treat the integral with an eighth-power mean-value estimate for each of the $S_j(\alpha)$,
\begin{align*}
\int_0^1 |S_j(\alpha)|^8 d\alpha & =  \int_0^1 \Big|\sum_{M_j^{1/3} < n \le N_j^{1/3}} \Lambda(n) e(a_j n^3 \alpha)\Big|^8 d\alpha \\
& =  \sum_{\substack{M_j^{1/3} < n_i \le N_j^{1/3} \\ i=1,\dots, 8}} \Lambda(n_1)\cdots \Lambda(n_8) \\  & \ \ \ \ \times  \int_0^1 e(a_j(n_1^3+n_2^3 + n_3^3 + n_4^3 - n_5^3 - n_6^3 - n_7^3 - n_8^3)\alpha)d\alpha\\ 
& =  \sum_{M_j^{1/3} < n \le N_j^{1/3}} [\Lambda(n)]^8 \sum_{\substack{ n_1^3 + n_2^3 + n_3^3 +  n_4^3 \\ = n_5^3 + n_6^3 + n_7^3 +  n_8^3 \\
M_j^{1/3} < n_k \le N_j^{1/3},\ k = 1,\dots, 8}} 1 \\
& \ll  \sum_{M_j^{1/3} < n \le N_j^{1/3}} \Lambda^8(n)  N_j^{4/3} \\
& \ll N_j^{4/3}L^7 \sum_{n \le N_j^{1/3}} \Lambda(n) \\
& \ll N_j^{5/3}L^7,
\end{align*}
where the last bound is from the Prime Number Theorem. 

In combination with H\"older's inequality gives
\begin{eqnarray} \label{eq:minor1}
\int_0^1 |S_1(\alpha)\cdots S_8(\alpha)|d\alpha & \le & \prod_{j=1}^8 \left(\int_0^1 |S_j(\alpha)|^8 d\alpha\right)^{1/8} \nonumber \\
& \ll & \frac{N^{5/3} L^7}{|a_1\cdots a_8|^{5/24}}.
\end{eqnarray}
For the remaining bound for $S_9(\alpha)$, we have from summation by parts,
\begin{eqnarray*}
\sum_{M_9^{1/3} < p \le N_9^{1/3}} e(a_9p^3 \alpha) & = & \int_{M_9^{1/3}}^{N_9^{1/3}} \frac{1}{\log t} dR_9(\alpha, t) \\
& = & \frac{R_9(\alpha, t)}{\log t} \bigg|_{M_9^{1/3}}^{N_9^{1/3}} + \int_{M_9^{1/3}}^{N_9^{1/3}} \frac{R_9(\alpha, t)}{t(\log t)^2}dt \\
& = & \frac{S_9(\alpha)}{\log N_9^{1/3}} + \int_{M_9^{1/3}}^{N_9^{1/3}} \frac{R_9(\alpha, t)}{t(\log t)^2}dt,
\end{eqnarray*}
where 
\[ R_9(\alpha, t) = \sum_{M_9^{1/3} < p \le t} (\log p) e(a_9 p^3 \alpha). \]
Hence, by Theorem 1 in Kumchev \cite{Kumchev1}, 
\begin{eqnarray*}
S_9(\alpha) & \ll & L(N_9^{(1/3)(1-1/14) + \epsilon} + N_9^{1/3+\epsilon}q^{-1/6}) \\
& \ll & L(N_9^{13/42+\epsilon} + N_9^{1/3+\epsilon}P^{-1/6}).
\end{eqnarray*}
By \eqref{eq:PQ}, we derive
\begin{eqnarray} \label{eq:minor2}
S_9(\alpha) & \ll & L(N_9^{13/42+\epsilon} + N_9^{1/3+\epsilon}N_j^{-1/60-\epsilon}) \nonumber \\
& \ll & LN_9^{19/60}.
\end{eqnarray}
Combining \eqref{eq:minor1} and \eqref{eq:minor2}, \eqref{eq:minor0} becomes
\begin{eqnarray} \label{eq:minorfinal}
\int_\MINA S_1(\alpha)\cdots  S_9(\alpha) e(-n\alpha ) d\alpha =  O\left(\frac{L^c N^{5/3 + 19/60}}{|a_1\cdots a_8|^{5/24} |a_9|^{19/60}}\right).
\end{eqnarray}

\section{Proofs of Theorems 3 and 4}
In view of (3.2), Lemma 3.2.4 and (3.63), we have
\[
r(n)=\frac{1}{3^9}\SS(n)\SI(n) + o\left( \frac{N^2}{|a_1\cdots a_9|^{1/3}}\right) +  O\left(\frac{L^c N^{5/3 + 19/60}}{|a_1\cdots a_8|^{5/24} |a_9|^{19/60}}\right).
\]
Comparing \eqref{eq:minorfinal} with the main term in \eqref{eq:SIorder}, if $n = N$ and all the $a_j$'s are positive, then 
\[\frac{N^2}{|a_1\cdots a_9|^{1/3}} \gg \frac{N^{5/3+19/60L^c}}{|a_1\cdots a_8|^{5/24}|a_9|^{19/60}}, \]
\begin{align*} {\text{i.e.,}} &\qquad N^{1/60} \gg |a_1\cdots a_8|^{1/8}|a_9|^{1/60}L^c \gg D^1D^{1/60+\epsilon}, \\
 {\text{i.e.,}} & \qquad N \gg D^{61+\epsilon}. \end{align*}
On the other hand, if not all of the $a_j$'s are the same sign, and $N \ge C|n|$, then 
\begin{align*} a_9p_9^3 &\le n - a_1p_1^3 - \cdots -a_8 p_8^3 \\
& \le n - 8N,\\
  {\text{or}} \qquad \qquad a_9 p_9^3 &\ll n + D^{61+\epsilon}.  \end{align*}
Therefore, without any loss of generality, for all $1 \le j \le 9$, we have
 \[ \qquad p_j \ll n^{1/3} + D^{20+\epsilon}. \]








%







\chapter{Remarks and Future Directions}
The obvious goal to Waring's problem and Waring-Goldbach's problem is to decrease the number of terms $s$ required to represent $n$. One approach would be to incorporate sieve methods with the circle method into the argument. When all the coefficients are all precisely one, Br\"udern in \cite{Brudern1} reduced the problem to only four terms but with one of the terms an almost prime. The other direction is to change the machinery of the circle method to eliminate the need for a minor arc estimate altogether (interested reader see Heath-Brown \cite{HB1}).

The limit of the estimate comes by taking the maximum (trivial bound) of seven of the terms and using Cauchy's inequality to bound the remaining two terms. By means of the large sieve we can estimate the two terms with an $L_2$ estimate. Heuristically, if we can find a way to find a $L_4$ or $L_8$-type estimate that with the combination of H\"older's inequality decrease the overall error term, we can increase the size of $P$ and improve the estimate.

 \end{document}